\documentclass[prl,aps,amsmath,amssymb]{revtex4-2}

\usepackage{times}

\usepackage{graphicx,amssymb,amsmath}
\usepackage{graphicx}
\allowdisplaybreaks

\newenvironment{sciabstract}{%
\begin{quote} \bf}
{\end{quote}}

\def\R{\mathbb R}
\def\P{\mathbb P}

\def\la{\left(}
\def\lb{\left\langle}
\def\rb{\right\rangle}
\def\ra{\right)}

\def\sech{{\rm \, sech}}

\begin{document}

\title{L\'evy flights as an emergent phenomenon in a spatially extended system}
\author{Chunxi Jiao$^{1}$}
\author{Georg A. Gottwald$^{2}$}
\email[]{georg.gottwald@sydney.edu.au}

\affiliation{$^{1}$Lehrstuhl f\"ur Angewandte Analysis, RWTH Aachen, 52062 Aachen, Germany}
\affiliation{$^{2}$School of Mathematics and Statistics, The University of Sydney, 2006 NSW, Australia}

\date{\today}

\maketitle












\begin{sciabstract}
Anomalous diffusion and L{\'e}vy flights, which are characterized by the occurrence of random discrete jumps of all scales, have been observed in a plethora of natural and engineered systems, ranging from the motion of molecules to climate signals. Mathematicians have recently unveiled mechanisms to generate anomalous diffusion, both stochastically and deterministically. However, there exists to the best of our knowledge no explicit example of a spatially extended system which exhibits anomalous diffusion without being explicitly driven by L\'evy noise. 

We show here that the Landau-Lifshitz-Gilbert equation, a stochastic partial differential equation (SPDE), despite only driven by Gaussian white noise, exhibits superdiffusive behaviour. The anomalous diffusion is an entirely emergent behaviour and manifests itself in jumps in the location of its travelling front solution. Using a collective coordinate approach we reduce the SPDE to a set of stochastic differential equations (SDEs) driven by Gaussian white noise. This allows us to identify the mechanism giving rise to the anomalous diffusion as random widening events of the front interface. 
\end{sciabstract}


\section{Introduction}
Brownian diffusion is ubiquitous in nature and social systems. Its ubiquity can be understood by the central limit theorem which states that appropriately scaled sums of random variables converge in distribution to a normal distribution. This is beautifully illustrated by Robert Brown's observation in 1827 of the emerging erratic behaviour of small particles suspended on the surface of a drop of water: the heavy particle experiences the sum of many (almost) uncorrelated kicks by the much smaller and lighter water molecules, leading to the apparent erratic motion of the particle. However, numerous examples in nature and in social systems defy this simple diffusive behaviour and exhibit anomalous diffusion. In systems with anomalous diffusion the standard central limit theorem breaks down but can be replaced by a generalized central limit theorem such that appropriately scaled sums converge in distribution to so called $\alpha$-stable processes \cite{Applebaum,SamorodnitskyTaqqu}. A particular class of anomalous diffusion is superdiffusion or L\'evy flights in which diffusion is enhanced compared to normal Gaussian diffusion by the occurrence of discontinuous jumps of arbitrary size. The presence of such jumps implies a fat tail in the associated probability distribution which may strongly deviate from the normal Gaussian distribution \cite{ZaburdaevEtAl15}. Superdiffusion has been observed, for example, in the motion of metal clusters and large molecules across crystalline surfaces \cite{SanchoEtAl04}, in active intracellular transport \cite{CaspiEtAl00}, conformal changes in proteins \cite{ReuveniEtAl10}, in living cell migration \cite{DieterichEtAl08,WeigelEtAl11,MalmiKakkadaEtAl18}, in optically active media \cite{DouglassEtAl12}, active turbulence \cite{MukherjeeEtAl21}, in tracer diffusion in active suspensions \cite{KanazawaEtAl20}, in atmospheric particle dispersion \cite{Richardson26}, in ice core data \cite{Ditlevsen99}, finance ratios \cite{Podobnik08}, power grids \cite{SchaeferEtAl18} and in foraging strategies of animals \cite{ViswanathanEtAl96,ViswanathanEtAl99,BartumeusEtAl03,Fernandez03,EdwardsEtAl07,GetzSaltz08,SimsEtAl08,BartumeusLevin08,HumphriesEtAl12,SimsEtAl19} to name but a few.  In these examples the emergence of superdiffusion is often dependent on geometric constraints which favour jumps, see also \cite{CherstvyEtAl13,StellaEtAl23b,CecconiEtAl23}. In recent years mathematicians have made substantial progress in identifying dynamical mechanisms which may give rise to superdiffusion without any superimposed geometric or topological constraints. There are two known dynamical mechanisms for which a process defined by appropriately scaled sums of observables converges to an $\alpha$-stable process: either for regular observables of weakly chaotic intermittent dynamics with a sufficiently slow decay of correlations \cite{Zweimuller03,Gouezel04,TyranKaminska10a,TyranKaminska10b,MelbourneZweimueller12,ChevyrevEtAl20}, or alternatively, in the case of strongly chaotic dynamics with rapid decay of correlations, for unbounded observables \cite{Gouezel08}. The same dynamical ingredients also hold for stochastic systems \cite{KuskeKeller01,PenlandSardeshmukh12,ThompsonEtAl17,AiminoEtAl22}. For readers who may not be familiar with these dynamic mechanisms to generate anomalous diffusion, we provide an illustrative introduction in the Appendix. Whereas these mechanisms have been verified in numerous examples for discrete maps as well as for ordinary differential equations \cite{GottwaldMelbourne13c,JungZhang18,JungEtAl19,JungEtAl20,MelbourneVarandas20,GottwaldMelbourne21,Gottwald21}, there are no explicit examples of spatially extended systems exhibiting superdiffusion. It has been conjectured, however, that superdiffusion should in principle be observable in spatially extended systems. This conjecture is part of a universal perspective based on the ambient symmetry of the problem \cite{GottwaldMelbourne13,GottwaldMelbourne16}. Spatially extended systems with symmetry can under certain conditions be reduced to the dynamics along the symmetry group and the so called shape dynamics which is orthogonal to it. A simple example of this reduction is given by traveling wave solutions of PDEs with translational invariance in one spatial dimension: the group dynamics is given by a linear drift of the reference frame in physical space, and the shape dynamics is provided by an ordinary differential equation obtained by transforming into the co-moving frame of reference. Such a reduction leads to the following dynamical scenario: the shape variables evolve according to their own dynamics, independent of the group variables, and then drive the dynamics of the group variables. The key insight developed in \cite{GottwaldMelbourne13,GottwaldMelbourne16,ChevyrevEtAl20} is that the group variables can exhibit anomalous superdiffusion provided the shape dynamics evolves according to the dynamics outlined above. In the deterministic setting it was proven that this is indeed possible in anisotropic media and it was conjectured that it is possible in isotropic media with odd spatial dimension. Numerical results as well as some mathematical theory suggest that the decomposition into the group and shape dynamics can be carried over to dissipative SPDEs with symmetry \cite{CartwrightGottwald19,HamsterHupkes20}. It is argued that the noise can freely diffuse along the neutrally stable group orbits whereas it is confined in the contracting shape dynamics. However, no explicit examples of PDEs or SPDEs exhibiting superdiffusions were known. 

In this work we present the first explicit example of a spatially extended system which exhibits superdiffusion without explicitly being driven by $\alpha$-stable noise. The anomalous diffusion is an entirely emergent behaviour. In particular, we consider the well studied stochastic Landau-Lifshitz-Gilbert (sLLG) equation \cite{LandauLifshitz35,Gilbert55,Gilbert04} which describes magnetisation along a one-dimensional wire. We will show that the location of the magnetisation interface of the front solution experiences superdiffusion. We remark that this is different to the superdiffusion observed in spin transport in one-dimensional spin chains which is a discrete finite dimensional microscopic model underlying the sLLG equation. Such spin systems were shown to fall into the Kardar-Parisi-Zhang (KPZ) universality class \cite{Znidaric11,LjubotinaEtAl19,DeNardisEtAl20,Bulchandani20,BulchandaniEtAl20,BulchandaniEtAl21,WeiEtAl22}. In the finite-dimensional spin chains  superdiffusion is manifested in the hopping of spins across sites and is {\em{described}} in the KPZ equation through the particular scaling of the correlation function. Instead in this work superdiffusion is directly observable in the SPDE as a L\'evy flight of the location of the front solution. Similarly, in the context of fluid dynamics and hydrodynamic transport, PDEs were introduced to describe transport of solute particles in stratified porous media subject to random-shear which exhibit superdiffusion \cite{MatheronDeMarsily80,CecconiEtAl23}. The superdiffusive behaviour in these random-shear models, however, again does not manifests itself in the solution of the transport PDE. Instead, the PDE describes the probability density function of random particles subject to vertical shear which experience superdiffusive dispersion in the latidudinal direction. Moreover, the superdiffusion is pre-described by the particular choice of the vertical shear profile. In contrast, the sLLG equation we consider here does not comprise any built-in mechanism which may lead us to suspect any anomalous diffusion and the anomalous diffusion is in that sense an emergent property.

The paper is organized as follows: We present in Section~\ref{sec:model} the SPDE under consideration, discuss its underlying symmetries and front solutions. We further show results of numerical simulations of the SPDE illustrating that the location of the front solution experiences jumps. To allow for a more thorough investigation of the mechanism giving rise to anomalous diffusion in the sLLG equation we perform in Section~\ref{sec:cc} a reduction of the infinite dimensional sLLG equation to a finite-dimensional system of SDEs. Using long numerical solutions of the SDEs we perform several statistical tests to show that indeed the location of the front interface is an $\alpha$-stable process and establish the dynamical mechanism giving rise to the anomalous diffusive behaviour. We conclude with a discussion in Section~\ref{sec:discussion}. For completeness, we provide in the Appendix an illustration how statistical limit theorems can lead to anomalous diffusion and outline the two mechanisms which may give rise to superdiffusion.  



\section*{The stochastic Landau-Lifshitz-Gilbert equation}
\label{sec:model}
The stochastic Landau-Lifshitz-Gilbert equation on a finite interval $ [0,L] $
\begin{eqnarray}
	dm = -m \times \la  H(m) +\lambda m \times H(m) \ra dt + \sigma m \times g \circ dW, 
	\label{eq:sLLG}
\end{eqnarray}
with homogenous Neumann boundary condition and initial condition satisfying $ \Vert m(x,0) \Vert=1 $ for all $ x \in [0,L] $,
where
\begin{eqnarray*}
	H(m)=\frac{\partial^2m}{\partial x^2}-m_2e_2- m_3e_3 
\end{eqnarray*}
and where  $W$ is one-dimensional Brownian motion, describes the magnetisation vector $m(x,t) = (m_1(x,t),m_2(x,t),m_3(x,t))$ of a ferromagnetic nanowire of length $L$ subject to dissipation (parametrised by $\lambda>0$) and thermal fluctuations \cite{LandauLifshitz35,Gilbert55,Gilbert04}. The sLLG equation is norm preserving in the sense that for every $ T \in (0,\infty) $, its martingale solution $ (\Omega, \mathcal{F}, (\mathcal{F}_t)_{t \in [0,T]}, \P, W, m) $ satisfies $\Vert m(x,t)\Vert = 1$, a.e.- $(x,t)$, $ \P $-a.s. The deterministic LLG equation, i.e. (\ref{eq:sLLG}) with $\sigma=0$, supports stable stationary front solutions, which for $L\to \infty$ take the form  
\begin{equation}
	\begin{array}{l}
		m_{0}(x) = \begin{pmatrix} \tanh(x) \\ \sech(x) \\ 0  \end{pmatrix}
	\end{array},
	\label{eq:m0}
\end{equation}
and enjoys translational symmetry as well as rotational symmetry, in the sense that if $m(x,t)$ is a solution then so is $R_1[\psi]m(x-\phi,t)$ for any constant $\phi$ and with $R_1[\psi]$ describing rotation around the $m_1$-axis by an angle $\psi$. Whereas the inclusion of the noise term in (\ref{eq:sLLG}) preserves the translational symmetry it breaks the rotational symmetry unless $g=(1,0,0)$. For small values of $\sigma$, however, the dynamics remains nearly rotationally invariant for some time. 

For the numerical simulation of the sLLG equation (\ref{eq:sLLG}) we employ a Crank-Nicolson scheme where nonlinearities are treated in an Adams-Bashforth scheme which ensures the preservation of the modulus $\Vert m \Vert^2$ \cite{SerpicoEtAl01,Banas05,dAquinoEtAl05}; for the noise we employ an Euler-Maruyama scheme \cite{Lord}.

\begin{figure}
	\centering
	\includegraphics[width=0.5\linewidth]{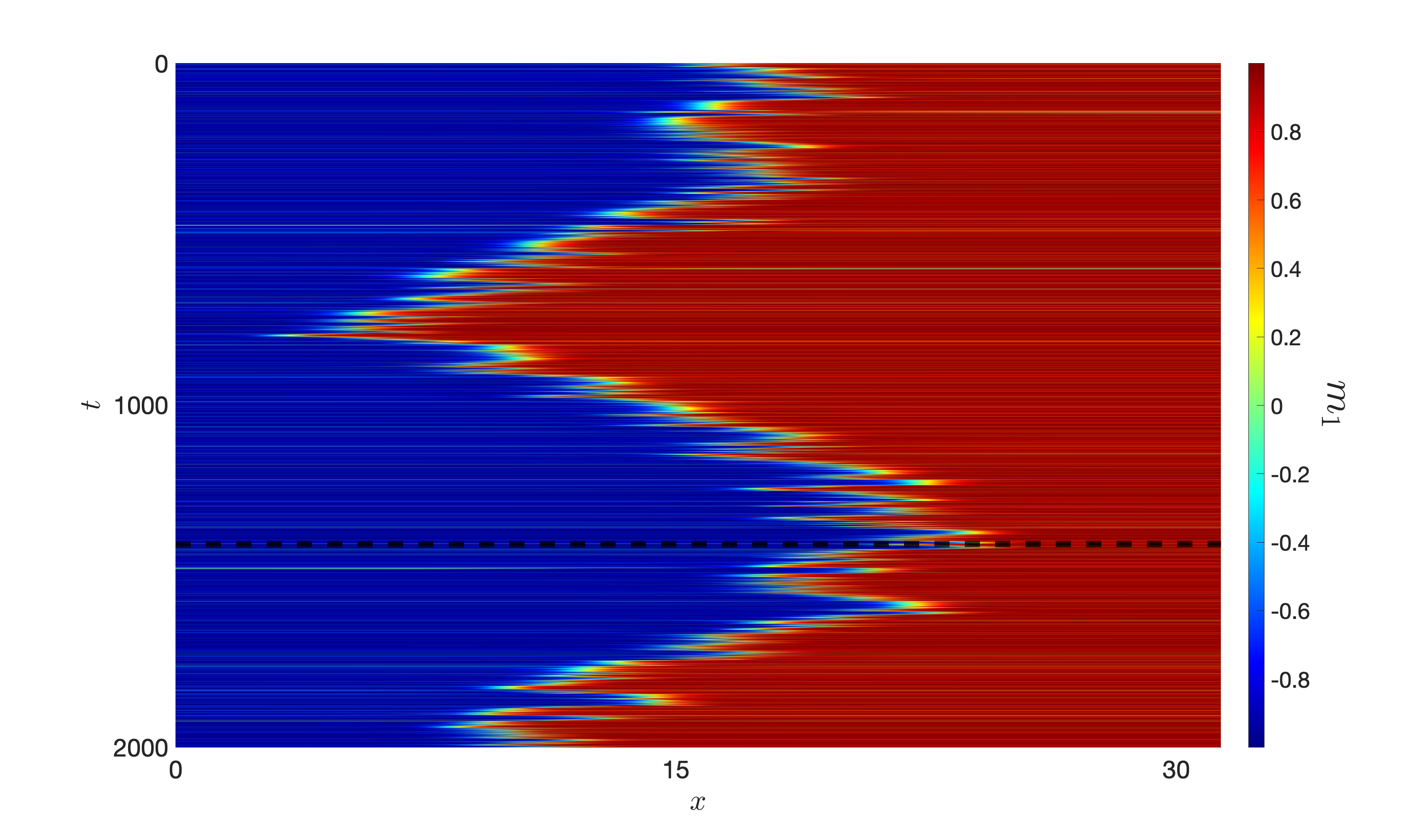}\\		
	\includegraphics[width=0.5\linewidth]{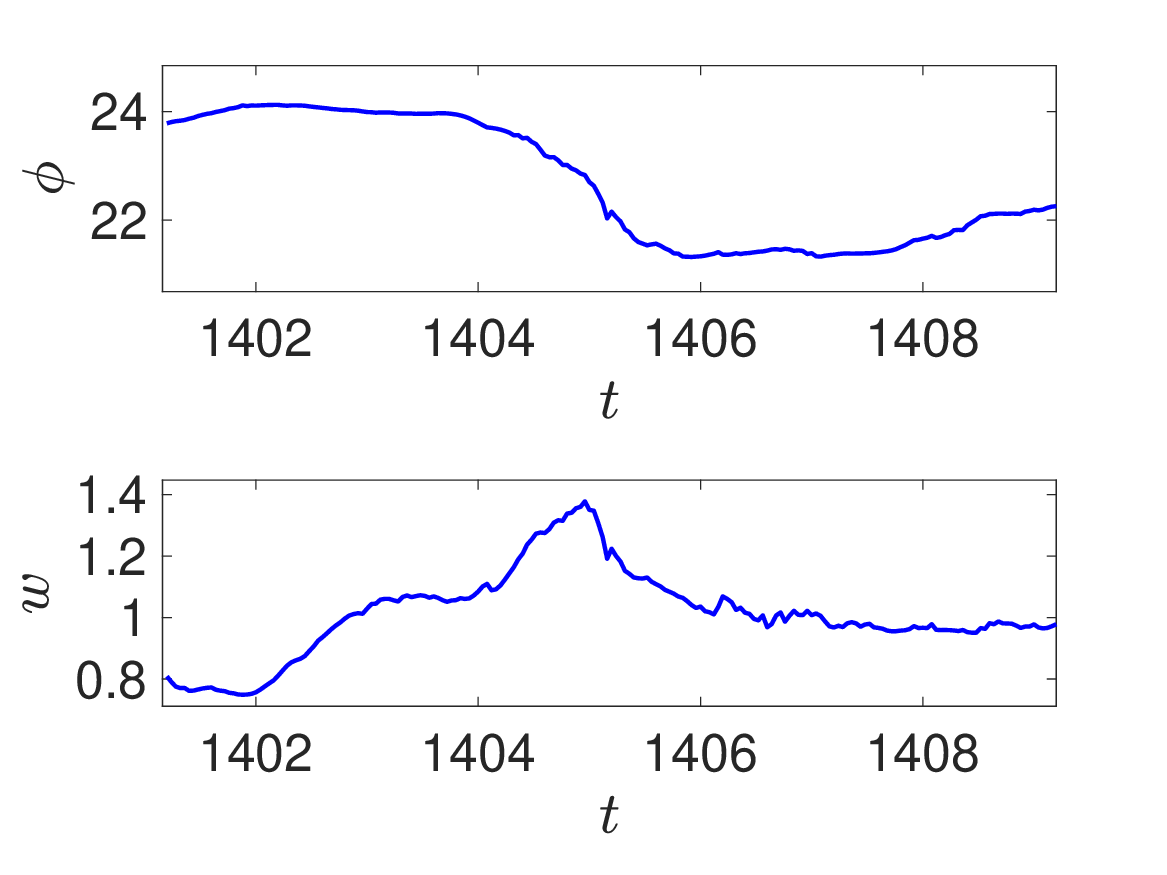}\\		
	\caption{Top: Contour plot of the solution $m_1$ of the sLLG equation (\ref{eq:sLLG}) for $\sigma=0.3$ and $\lambda=1$ with $g=(1,1,1)$ showing jumps in the location of the front interface. Zoom into a jump event at 
		$1404\lessapprox t \lessapprox 1406$ (marked as dashed horizontal lines in the contour plot) 
		showing the location of the front interface $\phi$ (middle) and the interface width $w$ (bottom).}
	\label{fig:m1}
\end{figure}

Figure~\ref{fig:m1} shows $m_1(x,t)$ for $\lambda=1$ and $\sigma=0.3$ with $g=(1,1,1)$. We further provide a zoom into a small time interval in which we show the temporal evolution of the location  $\phi(t)$ of the front interface and the width $w(t)$ of the front solution, which were both estimated by nonlinear least-square fitting (cf. (\ref{eq:ansatz})). It is seen that the front experiences what might be identified as a linear ballistic drift of its location $\phi(t)$ (albeit of only a small displacement of roughly 2 spatial units). These largish displacements of the front in a short time are accompanied by a sudden increase in the width $w(t)$ of the front solution. We will establish below that such events are indeed representative of anomalous diffusion and that the front interface location $\phi(t)$ experiences $\alpha$-stable L\'evy flights. It turns out, as we will show, that the length of the time window and size of the spatial domain necessary to observe significant jumps in simulations of the sLLG equation (\ref{eq:sLLG}), which would clearly identify the dynamics of the front interface as a L\'evy flight, are too large to be computationally feasible. To properly resolve the anomalous diffusion and analyse the superdiffusive dynamics we therefore now reduce the infinite dimensional sLLG equation to a set of finite-dimensional SDEs.



\section*{Reduction to an SDE via collective coordinates}
\label{sec:cc}

We shall employ the method of collective coordinates. Originally developed for deterministic PDEs  \cite{GottwaldKramer04} it has since been successfully employed to reduce the complexity of a wide range of dynamical systems including complex network dynamics \cite{Gottwald15,Gottwald17,SmithGottwald20} and SPDEs \cite{CartwrightGottwald19,CartwrightGottwald21}. In the context of SPDEs with symmetries the method of collective coordinates relies on a decomposition of the dynamics into the dynamics along the group and the dynamics orthogonal to it, the so-called shape dynamics. Assuming that the deterministic core supports a travelling front, as is the case for the sLLG equation, the question is how the noise will affect the solution. It is argued that the noise can freely diffuse along the neutrally stable group orbit, given here by the translational symmetry, whereas the noise is controlled in the contracting ($\lambda>0$) shape dynamics \cite{CartwrightGottwald19}. To capture this, we consider the ansatz 
\begin{equation}
	\label{eq:ansatz}
	\hat m(x,t) = R_3[\theta(t)]\,R_2[\eta(t)]\,R_1[\psi(t)]\, m_0(w^{-1}(t)(x-\phi(t)).
\end{equation}
For simplicity we assume from now on that $L\to \infty$ which allows to use the explicit solution (3.2) for $m_0$. 
Here $R_j[\xi]$ denotes rotations around the $j^{\rm{th}}$ axis by the angle $\xi$. The dynamics of the infinite-dimensional SPDE is now encapsulated in the time evolution of the shape variables $(w,\theta,\eta,\psi)$ and the translational group variable $\phi(t)$. Indeed, we are able to fit such an ansatz to solutions of numerical simulations with a high degree of accuracy (cf. Section \ref{Section: CC vs SPDE}). The rotations imply that each component of $\hat m(x,t)$ can be written as a linear combination of $\tanh$ and $\sech$ functions with the coefficients depending on the angle variables. These amplitudes evolve stochastically on the sphere defined by $\Vert \hat m \Vert^2=1$ (cf. Figure~\ref{fig:w}f).\\ 

The framework of collective coordinates views the ansatz (\ref{eq:ansatz}) as a Galerkin approximation and appropriate projection leads to an SDE for the collective coordinates. For expository ease we denote by $\vec{p}=(w, \theta, \eta, \phi, \psi)$ the vector of collective coordinates which are assumed to evolve according to a stochastic differential equation 
\begin{equation*}
	dp_i = f_{p_i}(\vec{p}) \ dt + \sigma_{p_i}(\vec{p}) \ dW, 
\quad  
i = 1, \ldots, 5. 
\end{equation*}
Our aim in the following is to determine the drift coefficients $f_{p_i}$ and the diffusion coefficients $\sigma_{p_i}$  \cite{CartwrightGottwald19}. Applying It{\^o}'s lemma to the ansatz (\ref{eq:ansatz}) yields 
\begin{align*}
	d\hat m &= \left[ \sum_{i=1}^5 f_{p_i} \partial_{p_i} \hat m + \frac{1}{2} \sum_{i,j=1}^5 \sigma_{p_i} \sigma_{p_j} \partial_{p_i, p_j} \hat m \right] dt + \sum_{i=1}^5 \sigma_{p_i} \partial_{p_i} \hat m \ dW. 
\end{align*}
Introducing the drift and diffusion terms of the sLLG equation \eqref{eq:sLLG}
\begin{align*}
	b_f(m) &= - m \times \la H(m) + \lambda m \times H(m) \ra + \frac{1}{2}\sigma^2(m \times g) \times g, \\
	b_\sigma(m) &= \sigma m \times g,
\end{align*}
we define the error process $ \{\mathcal{E}(t) : t \geq 0 \} $ as the residual  
\begin{align*}
	d\mathcal{E} &= \left[ \sum_{i=1}^5 f_{p_i} \partial_{p_i} \hat m + \frac{1}{2} \sum_{i,j=1}^5 \sigma_{p_i} \sigma_{p_j} \partial_{p_i p_j} \hat m -b_f(\hat m)  \right] dt 
	+ \left[\sum_{i=1}^5 \sigma_{p_i} \partial_{p_i} \hat m -b_\sigma(\hat m) \right] dW
\end{align*}
with $ \mathcal{E}(0)=0 $. If the ansatz (\ref{eq:ansatz}) were a solution of the sLLG equation, then $ d\mathcal{E} = d \mathcal{E}(x,t;\vec{p}) =0 $ for all $(x,t)$. 

In order to maximise the degree to which the ansatz $ \hat m = \hat m(\vec{p}) $ approximates the solution of \eqref{eq:sLLG}, we minimise the projection of the error $ \mathcal{E} $ onto the space spanned by the collective coordinates $\vec{p}$. We do so by requiring that the error is orthogonal to the tangent space of the ansatz manifold parametrised by the collective coordinates. In particular, $f_p$ and $\sigma_p$ are determined such that the error is orthogonal to the tangent space spanned by $ \{\partial_{p_k} \hat m \}_{k=1}^5 $, i.e. 
\begin{align*}
	\lb d\mathcal{E}, \partial_{p_k} \hat m \rb=0, 
	\quad  
k = 1, \ldots, 5, 
\end{align*}
where the angular brackets denote the inner product on $L^2([0,L];\R^3)$ with  
$\lb a, b \rb = \int_{0}^{L} a(x)b(x)dx$. Separating into drift terms which are proportional to $dt$ and diffusion terms which are proportional to $dW$ we require
\begin{align}
	\label{eq:Ik}
	\lb \sum_{i=1}^5 \sigma_{p_i} \partial_{p_i} \hat m , \partial_{p_k} \hat m \rb = \lb b_\sigma(\hat m) , \partial_{p_k} \hat m  \rb=0, 
\end{align}
and
\begin{align}
	\label{eq:Hk}
	&\lb \sum_{i=1}^5 f_{p_i} \partial_{p_i} \hat m + \frac{1}{2} \sum_{i,j=1}^5 \sigma_{p_i} \sigma_{p_j} \partial_{p_i p_j} \hat m ,\partial_{p_k} \hat m \rb  
	= \lb b_f(\hat m) ,\partial_{p_k} \hat m \rb=0 
\end{align}
for each $ k=1,\ldots,5 $. This constitutes a linear system of 10 equations for the 10 unknown drift and diffusion terms $f_{p_i}$ and $\sigma_{p_i}$. In practice, the diffusion coefficients $\sigma_{p_i}$ are determined by \eqref{eq:Ik}, and then inserted into \eqref{eq:Hk} to yield a linear system for the drift terms $f_{p_i}$. 
We first solve the system for finite domain length $L$ and then take the limit $L\to \infty$. Using the software package Mathematica \cite{Mathematica}, taking $L\to \infty$ leads 
%
to the following explicit expressions for the dynamics of the collective coordinates $(\psi,\theta,\eta,w)$ and $\phi$
\begin{eqnarray}
	{\text{(shape dynamics) }}
	\begin{cases}
		dw = f_w(\theta,\eta,\psi,w)\, dt\\
		d\theta = f_\theta(\theta,\eta)\, dt + \sigma_\theta dW_t\\
		d\eta = f_\eta(\theta,\eta)\, dt + \sigma_\eta dW_t\\
		d\psi = f_\psi(\theta,\eta,\psi,w)\, dt + \sigma_\psi dW_t
	\end{cases}
	\label{eq:cc_s}
\end{eqnarray}
\begin{eqnarray}
	{\text{(group dynamics) }}
	\begin{cases}
		d\phi = f_\phi(\theta,\eta,\psi,w)\, dt, 
		\hphantom{\qquad\qquad}
         \end{cases}
	\label{eq:cc_g}
\end{eqnarray}
with the drift and diffusion terms 
\begin{eqnarray*}
	f_w &=& \frac{6}{\pi^2} \lambda  w^{-1} - \frac{3}{2\pi^2} w \Big[\lambda  \left(\cos^2\theta \cos 2\eta (3-\cos 2\psi) -2 \sin 2\psi \sin 2\theta \sin \eta \right. \nonumber \\
	&&\left. +\cos^2\psi (3 \cos 2\theta-1)\right)+2 \pi  \cos \eta \left(2 \cos \psi \cos^2\theta \sin \eta-\sin \psi \sin 2\theta\right) \Big]
	\label{eq:SMCC1}\\
	f_\theta &=& \frac{1}{2} \left(\sigma^2 (\tan \eta  (2 \cos 2 \theta  \tan \eta +\sin \theta -\cos \theta )+\cos 2 \theta )-\lambda  \sin 2 \theta \right) \\
	f_\eta &=& \lambda  \cos^2\theta  \sin \eta  \cos \eta +\frac{1}{2} \sigma^2 (\sin \theta +\cos \theta ) (\tan \eta  (\sin \theta +\cos \theta )-1) \\
	f_\phi &=& \frac{w}{2} \la \cos \theta \sin \psi \sin \eta + \cos \psi \sin \theta  \ra \big( 2 \cos \psi \cos \theta \sin \eta - 2 \sin \psi \sin \theta \\
	&&\quad + \lambda \pi \cos \theta \cos\eta \big) \\
	f_\psi &=& \frac{\sigma^2}{2} \sec \eta \la \cos \theta - \sin \theta - 2 \cos 2\theta \tan \eta \ra \\
	&&\quad + \frac{\lambda}{8} \la 3 \cos 2\theta - 2 \cos^2 \theta \cos 2 \eta -1\ra \sin 2\psi + \lambda \cos^2 \psi \sin 2\theta \sin \eta
	\label{eq:SMCC2}
\end{eqnarray*}
and the diffusion terms are 
\begin{eqnarray*}
	\sigma_\theta &=& \sigma (\tan \eta  (\sin \theta +\cos \theta )-1) \label{eq:SMCC1b}\\
	\sigma_\eta &=& \sigma (\sin \theta -\cos \theta ) \\
	\sigma_\psi &=& -\sigma \sec \eta (\cos \theta + \sin \theta).
	\label{eq:SMCC2b}
\end{eqnarray*}

\subsection{Ability of the collective coordinate approach to reproduce the dynamics of the full SPDE}\label{Section: CC vs SPDE}
We provide some evidence that the collective coordinate framework is indeed able to reproduce the dynamics of the full sLLG equation (\ref{eq:sLLG}). We first show that the collective coordinate approach is able to capture the statistical behaviour of the full SPDE. 

\begin{figure}[h]
	\centering
	\includegraphics[width=0.32\linewidth]{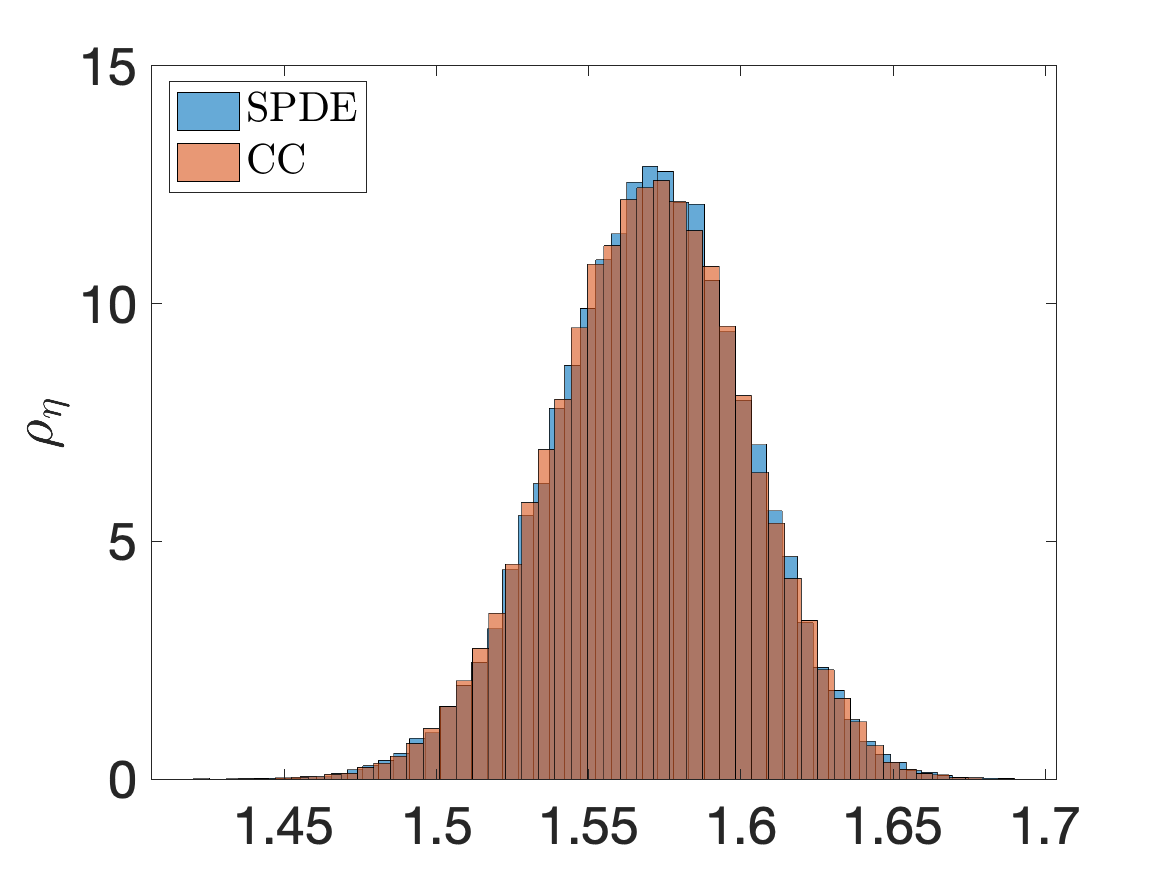}\;
	\includegraphics[width=0.32\linewidth]{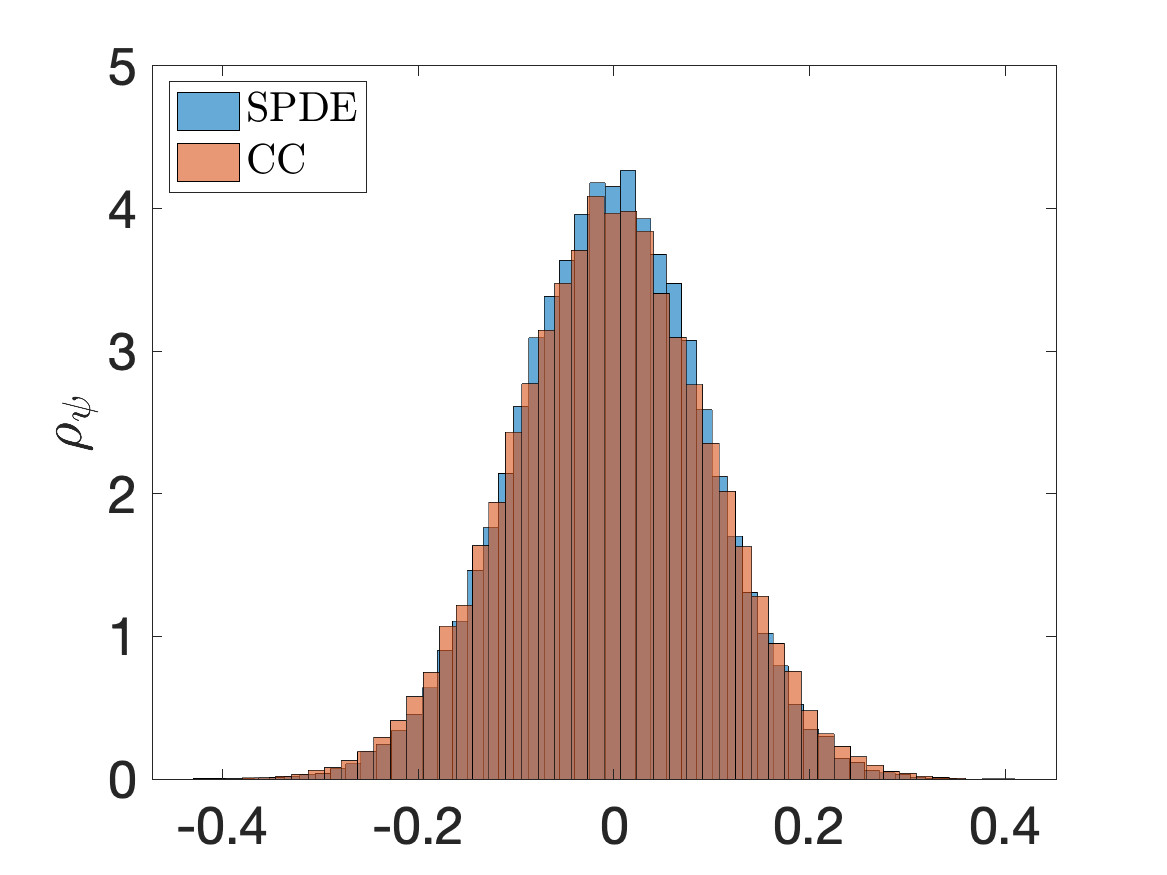}\;
	\includegraphics[width=0.32\linewidth]{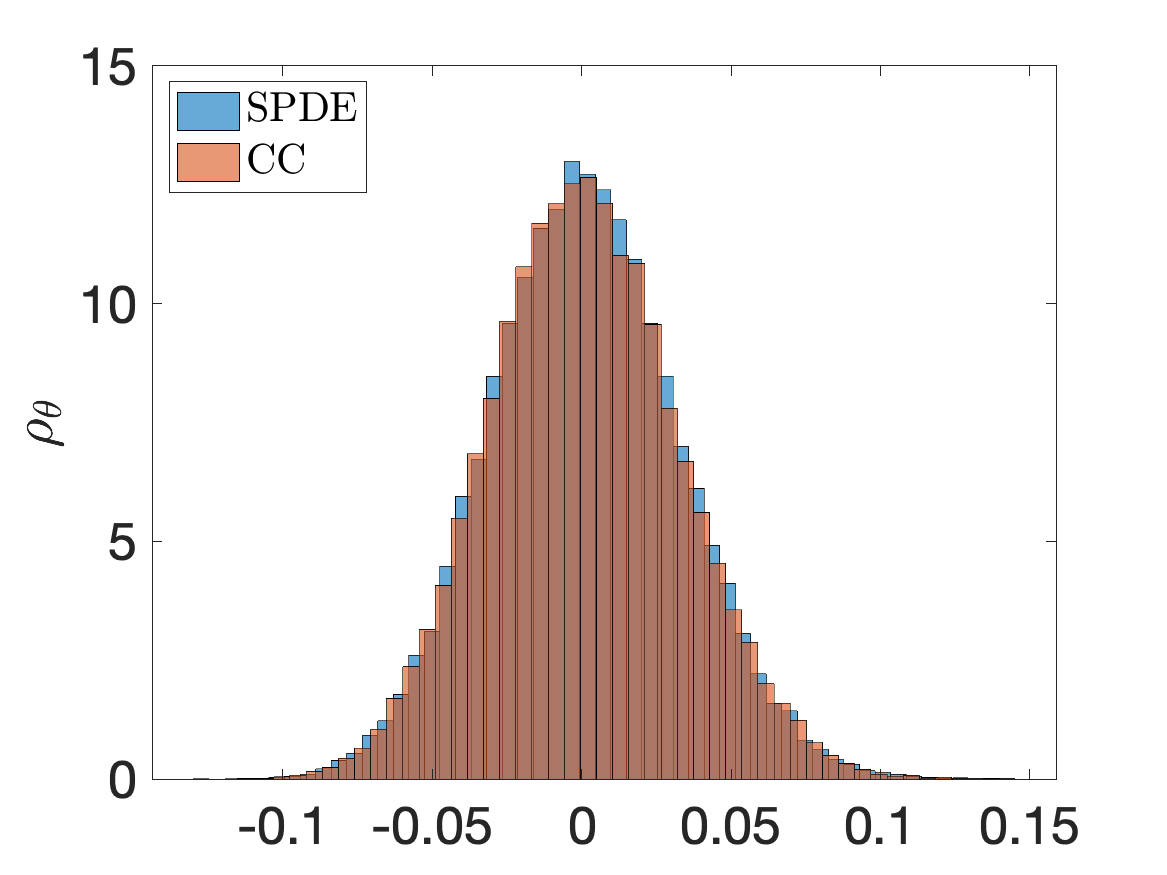}
	\caption{Comparison of the normalized empirical histograms of the rotational variables between the full SPDE (\ref{eq:sLLG}) 
	and the reduced collective coordinate system fitted (\ref{eq:cc_s})-(\ref{eq:cc_g}) at $t=1$. Parameters are $\lambda=5$ and $\sigma=0.1$ for $g=(1,1,1)$. The values corresponding to simulations of the SPDE were obtained through a nonlinear least square fit.}
	\label{fig:pdf}
\end{figure}

Figure~\ref{fig:pdf} shows a comparison of the empirical probability function of the angle variables $(\eta,\theta,\psi)$ at fixed time $t=1$ for $\lambda=5$, $\sigma=0.1$ with excellent agreement between the collective coordinate solutions and the full SPDE (on the bounded domain $[-\frac{L}{2}, \frac{L}{2}]$ with $L=10 \pi$). The angle variables for the full SPDE were obtained using a nonlinear least-square fit to the ansatz for the front solution (\ref{eq:ansatz}).\\

To provide further indication about the reliability of the collective coordinate description of the dynamics of front interfaces in the sLLG equation we show in the following that the collective coordinate ansatz recovers an analytical solution in the special case of $g=(1,0,0)$. 

In the special case when $g=(1,0,0)$ the sLLG equation exactly preserves the rotational symmetry along the axis of rotation aligned with $m_1$. In this case the sLLG equation (\ref{eq:sLLG}) on $\R$ with initial condition $ m(\cdot,0) = m_0 $ given by \eqref{eq:m0} permits the analytical solution $m(x,t) = e^{G W(t)} m_0(x)$. 
Indeed, for $g=(1,0,0)$ the multiplicative noise of the sLLG equation can be written as $\sigma m\times g \circ dW = G m \circ dW$ with the matrix  
\begin{align*}
	G = \sigma \begin{pmatrix} 0 & 0 & 0 \\ 0 & 0 & 1 \\ 0 & -1 & 0  \end{pmatrix}.
\end{align*}
Then 
\begin{align}
m(x,t) = e^{G W(t)} m_0(x) 
\label{eq:m_gG}
\end{align}
is an explicit solution of  (\ref{eq:sLLG}). This readily follows upon application of It{\^o}'s lemma, according to
\begin{align*}
	dm
	&= e^{GW} dm_0 + Ge^{GW} m_0 \circ dW \\
	&= e^{GW} (-m_0 \times H(m_0) - \lambda m_0 \times (m_0 \times H(m_0))) \ dt + G m \circ dW \\
	&= -\la e^{GW} m_0 \times H(e^{GW} m_0) + \lambda e^{GW} m_0 \times \la e^{GW} m_0 \times H(e^{GW} m_0) \ra \ra dt 
	+ Gm \circ dW \\
	&= -\la m \times H(m) + \lambda m \times \la m \times H(m) \ra \ra \ dt + \sigma m \times g \circ dW(t),
\end{align*} 
which is the sLLG equation (\ref{eq:sLLG}). Here the commutation of $G$ and $H(m)$ follows from the particular form of $G$ for $g=(1,0,0)$. In particular, $e^{GW(t)} m_0 = R_1[-\sigma W(t)] m_0$, i.e. the dynamics is a random motion around the $m_1$-axis. This implies that $m_1(x,t) = \tanh(x)$ and that $m_{2,3}(x,t) = A_{2,3}(t)\sech(x)$ with $A_{2,3}(t)$ being stochastic processes with the constraint that $\Vert m(x,t)\Vert^2=1$. 

The collective coordinate system (\ref{eq:cc_s})-(\ref{eq:cc_g}) recovers the exact solution (\ref{eq:m_gG}).  
In the case of rotational symmetry along the axis of rotation aligned with $m_1$, denoted by $R_1$, we set $\theta=\eta=0$ in our ansatz (\ref{eq:ansatz}), leading to the following equations for the collective coordinates
\begin{eqnarray}
	{\text{(shape dynamics) }}
	\qquad
	\begin{cases}
		dw = \frac{6\lambda}{\pi^2 w}(1-w^2) dt
	\end{cases}
	\label{eq:grot1}
\end{eqnarray}
\begin{eqnarray}
	{\text{(group dynamics) }}
	\qquad 
	\begin{cases}
		d\psi = -\sigma dW(t)\\
		d\phi = 0
	\end{cases}.
	\quad \;
	\label{eq:grot2}
\end{eqnarray}
The equation for the interface width converges to the stable fixed point $w=1$ and the group variables are readily solved to yield $\psi(t) = \psi(0)-\sigma W(t) $ and $ \phi(t) = \phi(0) $ for all $ t $. Thus under the initial condition $ (w(0), \psi(0), \phi(0)) = (1,0,0) $, the ansatz reduces to $ \hat m(x,t) = R_1[-\sigma W(t)] m_0(x) $, recovering the exact solution. Note that in this case with $g=(1,0,0)$ the rotational symmetry around the $m_1$-axis implies that $\psi$ is a group variable rather than a shape variable. Hence the collective coordinate system (\ref{eq:grot1})--(\ref{eq:grot2}) is of the form where the noise is only active in the neutrally stable group dynamics.


\subsection{Mechanism for anomalous diffusion}
\label{sec:mechanism}

The angles $(\psi,\theta,\eta)$ ensure that $\Vert m\Vert^2=1$, i.e. the solution remains on the sphere, and the dynamics can be viewed as the axis of this sphere performing a random walk. To identify the mechanism leading to superdiffusion, it is instructive to express the ansatz (\ref{eq:ansatz}) as 
\begin{align}
\label{eq:ansatzAB}
	m(x,t) &= A(t) \tanh(w^{-1}(t)(x-\phi(t))) + B(t)\sech(w^{-1}(t)(x-\phi(t))),
\end{align}
where 
\begin{align*}
	A = \begin{pmatrix} \cos \theta \cos \eta \\ \sin \theta \cos \eta \\ \sin \eta \end{pmatrix}
	,\; 
	B = \begin{pmatrix} -\cos \psi \sin \theta - \cos \theta \sin \psi \sin \eta \\ \cos \psi \cos \theta - \sin \psi \sin \theta \sin \eta \\ \cos \eta \sin \psi \end{pmatrix}.
\end{align*}
The evolution equations 
for the collective coordinates $(\psi,\theta,\eta,w,\phi)$  
can then be written in terms of the shape variables $(A,B,w)$ and the group variable $\phi$. For $g=(1,1,1)$ we obtain for the shape dynamics
\begin{align}
	dA &= \lambda A \times (A \times KA) dt + \sigma A \times g \circ dW(t), \label{eq:CCA}\\
	dB &= \lambda B \times \Big[ \langle B \times KB, A \rangle A  + \frac{A_1^2 A_3}{1-A_3^2} A \times e_3 
	+ \frac{A_1A_2}{1-A_3^2} (e_3 - A_3 A) \Big] dt \nonumber \\ & \quad 
	+ \sigma B \times g \circ dW(t),  \label{eq:CCB}\\
	dw &= \frac{6}{\pi ^2 w} \left[\pi  w^2 A_1 \lb A\times B, e_1 \rb - \lambda  w^2 \left(\lb A \times B, e_2 \rb^2-A_3^2-B_1^2+B_2^2\right)+\lambda \right] dt
	\label{eq:CCw}	
\end{align}
and for the group dynamics
\begin{eqnarray}
	d\phi = w B_1 \left[  \lb A\times B, e_1 \rb -\frac{\pi}{2} \lambda A_1 \right] dt ,
	\label{eq:CCphi}
\end{eqnarray}
where $K = {\rm{diag}{(0,1,1)}}$. 



\begin{figure}[h]
	\centering
	\includegraphics[width=0.5\linewidth]{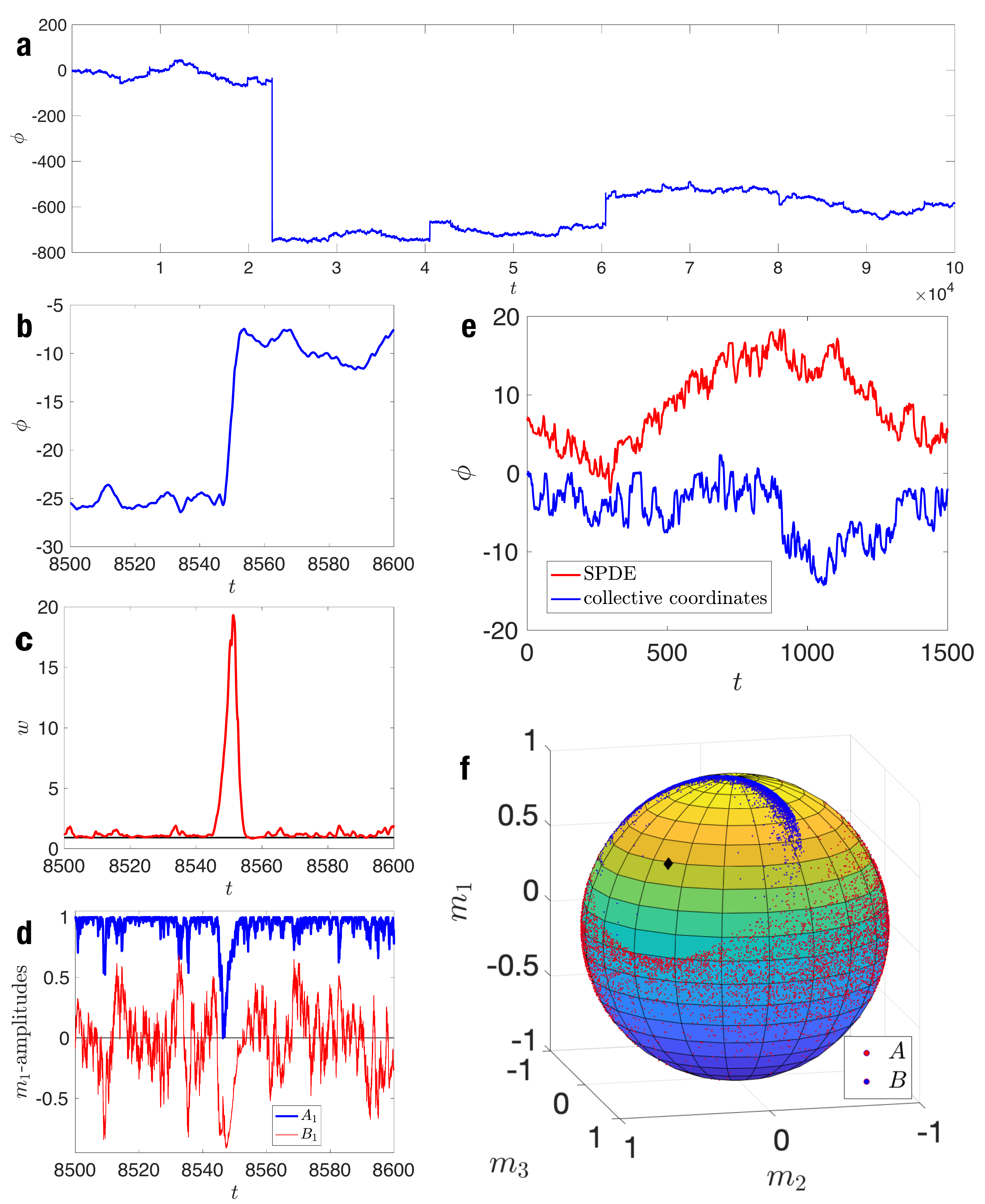}
	\caption{Solutions of the reduced collective coordinate system (\ref{eq:CCA})-(\ref{eq:CCphi}) for $\lambda=1$ and $\sigma=0.3$ with $g=(1,1,1)$. (a): Translational symmetry variable $\phi(t)$. (b)-(d): Mechanism generating jumps: Zoom into a single jump (b) of $\phi$, which is caused by an extreme interface width $w$ (c), which itself is generated by random events of near-zero amplitude of the $\tanh$-component $A_1$ of $m_1$ (d). (e): Typical time series of the front interface location $\phi(t)$ obtained from simulations of the sLLG (\ref{eq:sLLG}) and of the reduced collective coordinate system (\ref{eq:CCA})-(\ref{eq:CCphi}). We arbitrarily shifted the time series to allow for a comparison.  (f): Time series of the amplitudes $A_{1,2,3}$ and $B_{1,2,3}$ on the sphere with radius $1$. The black diamond marks where $g=(1,1,1)$ pierces the sphere.}
	\label{fig:w}
\end{figure}

\begin{figure}[h]
	\centering
	\includegraphics[width=0.5\linewidth]{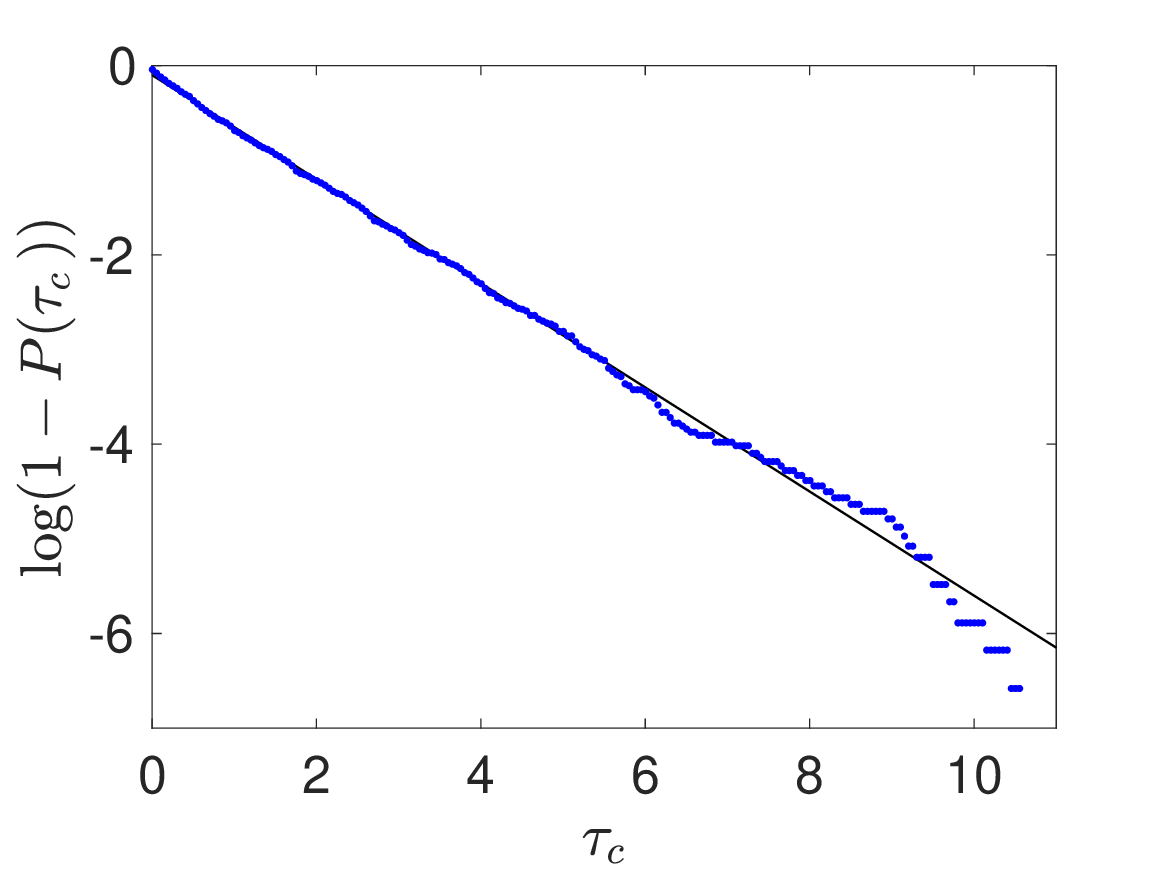}
	\caption{Log plot of the empirical normalised cumulative probability density function of waiting times $\tau_c$ between $w>1$ events (online blue). The continuous line with slope $-0.55$ depicts the cumulative probability density function of a Poisson process with exponentially distributed waiting time as a best fit. Parameters are $\lambda=1$ and $\sigma=0.3$ and $g=(1,1,1)$ as in Figure~\ref{fig:w}. }
	\label{fig:Poisson}
\end{figure}

The collective coordinate system (\ref{eq:CCA})-(\ref{eq:CCphi}) allows for long time integrations and enables us to study how anomalous diffusion is generated. In Figure~\ref{fig:w} we show results of a numerical simulation of (\ref{eq:CCA})-(\ref{eq:CCphi}) illustrating the chain of events leading to superdiffusion of the translational variable $\phi$ seen in Figure~\ref{fig:w}(a). 
Random incidences of (near) flips of the magnetization front $m_1$ with $A_1\approx 0$ and $B_1\approx \pm 1$, as seen in Figure~\ref{fig:w}(d) \footnote{This corresponds to $\theta=\tfrac{\pi}{2}$ (or $\theta=\tfrac{3\pi}{2}$) and $\psi=0$ (or $\psi=\pi$).}, lead to distinct spikes in the front interface width $w$ (see Figure~\ref{fig:w}(c); cf. Figure~\ref{fig:m1} for the full SPDE over a shorter time period). 
Events of $w$ exceeding specified thresholds are Poisson distributed, suggesting that large $w$-events and the generating (near) flips are uncorrelated. Indeed, Figure~\ref{fig:Poisson} shows that the cumulative probability distribution function of the waiting times $\tau_c$ between events of $w\ge 1$ closely resembles that of a Poisson process (with exponential inter-arrival times) where $P(\tau_c) = 1 - \exp(-\tau_c/\bar{\tau_c})$ with $\bar{\tau_c} = 0.55$.
These (arbitrarily large) spikes in $w$ are then causing the jumps in the interface location $\phi$ generating a L\'evy or $\alpha$-stable process (see Figure~\ref{fig:w}(b)). This chain of events is clearly reflected in equations (\ref{eq:CCA})-(\ref{eq:CCphi}): for $A_1\approx 0$ and $B_1\approx \pm 1$ a Taylor expansion yields $dw=6\lambda (\pi^2 w)^{-1}(1+(1-B_2^2)w^2)\, dt>0$. This implies an increase in $w$ which in turn leads to $d\phi \sim \pm w \, dt$ causing large jumps in $\phi$, the signs of which being determined by the random sign of $B_1[ \lb A\times B, e_1 \rb - \pi \lambda A_1/2]$. 

The anomalous diffusion in the translational group variable $\phi$ is generated by the shape dynamics. However, none of the driving shape dynamics of $(w,\theta,\eta,\psi)$ exhibit fat tails. The anomalous diffusion is instead caused by $d\phi$ being unbounded for $w\to \infty$. Hence using the framework of collective coordinates we identified that in the sLLG equation superdiffusion is generated by the known mechanism of unbounded driving observables, illustrated in the Appendix. 
Figure~\ref{fig:w}(e) shows that on shorter time scales $\phi$ obtained by the collective coordinate ansatz closely resembles $\phi$ obtained from simulating the full SPDE.

\subsection{Testing for superdiffusion}
\label{sec:test}
L\'evy processes or $\alpha$-stable processes are characterized by three parameters $(\alpha,\beta,\gamma)$. The parameter $0<\alpha\le 2$ controls the anomalous character of the diffusion with $\alpha=2$ being normal diffusion. The skewness parameter $-1\le \beta \le 1$ describes the probability of the process experiencing up- or downward jumps with $\beta=\pm1$ allowing only for up-/downward jumps, respectively. The parameter $\gamma$ is a scale parameter. 

\begin{figure}[h]
	\centering
	\includegraphics[width=0.5\linewidth]{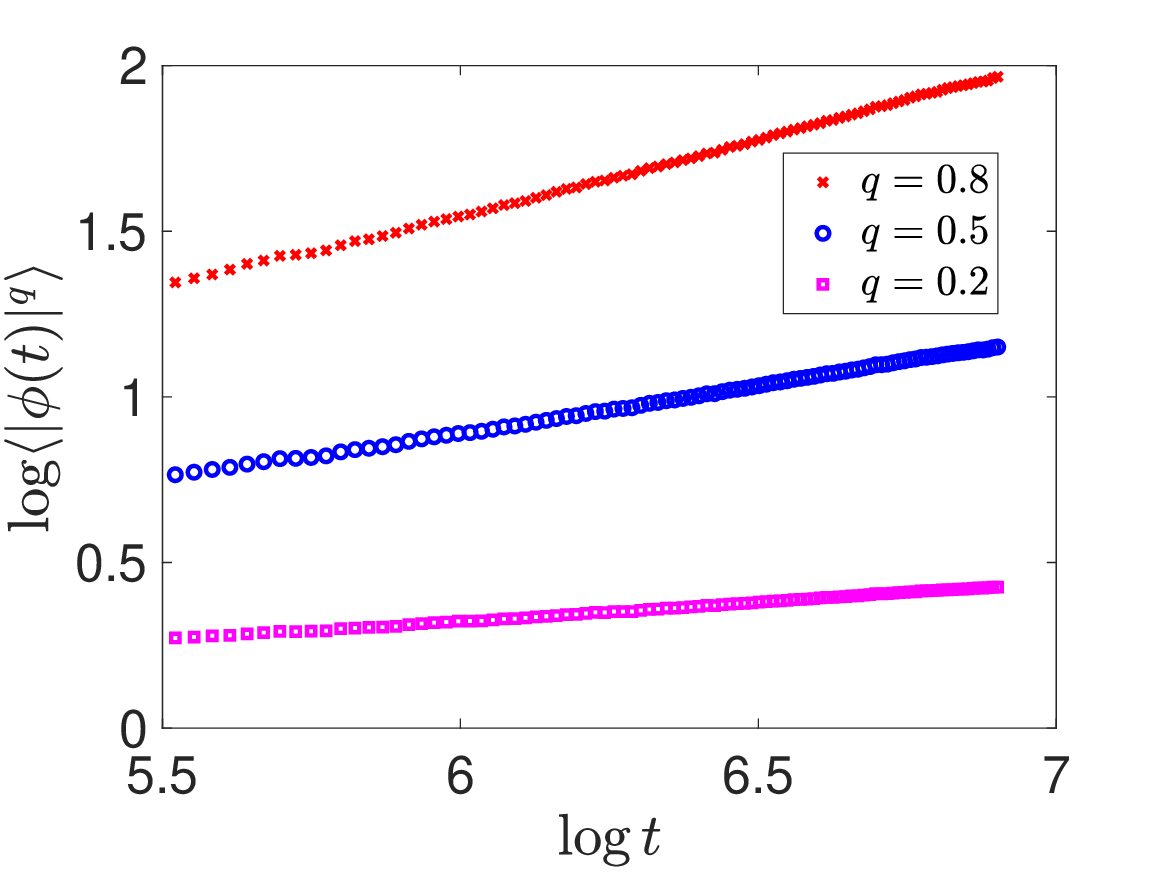}
	\caption{Log-log plot of the $q$th moment $\langle |\phi(t)|^q\rangle$ versus time $t$ for three values of $q$. For $q=0.8$, $q=0.5$ and $q=0.2$ the best linear fit has slopes $\kappa=0.4592$, $\kappa=0.2873$ and $\kappa=0.1151$, respectively, all suggesting $\alpha=q/\kappa=1.74$ (cf. (\ref{eq:MSQ})). To compute the ensemble average $10,000$ realizations were used.}
	\label{fig:MSQ}
\end{figure}

To test if the dynamics of the translational variable $\phi(t)$ is driven by an $\alpha$-stable process and to determine the value of $\alpha$ we employ several test statistics. We begin by determining the scaling behaviour of the $q$th moment $\langle |\phi(t)|^q\rangle $ with time (note that $\phi$ has zero mean). The $q$th moment is defined for $\alpha<q$ and, for stochastic dynamics, scales like 
\begin{align}
\langle |\phi(t)|^q\rangle \sim t^\frac{q}{\alpha},
\label{eq:MSQ}
\end{align}
where the angular bracket denotes an ensemble average \cite{SamorodnitskyTaqqu}. For anomalous diffusion with $\alpha<2$ the second moment is not defined and an ensemble average does not converge with an increasing number of ensemble members. Figure~\ref{fig:MSQ} shows the scaling of the $q$th moment $\langle |\phi(t)|^q\rangle$ in a $\log$-$\log$ plot for $q=0.2$, $q=0.5$ and $q=0.8$ with a very clear linear scaling. The slope of a best linear fit suggests $\alpha=1.74$ for all three values of $p$ indicative of superdiffusion. We remark that the scaling of the $q$th moment is typically a better test statistics for smaller values of $q$ (see, for example, \cite{GottwaldMelbourne16b}).

To further test for $\alpha$-stable behaviour, we additionally employ the $p$-variation \cite{MagdziarzEtAl09,MagdziarzKlafter10,HeinImkellerPvalyukevich}, a test for superdiffusion based on the scaling behaviour of low-order moments \cite{GottwaldMelbourne16b}, the maximum likelihood test \cite{DuMouchel73,Nolan01}, a quantile matching test \cite{Mcculloch86} and employing linear regression of the empirical characteristic function \cite{Koutrouvelis80,Koutrouvelis81}. Unlike simply testing for fat tails in the empirical histogram, the presence of which is not necessarily linked to an underlying $\alpha$-stable process, the above methods are proper statistics to identify $\alpha$-stable behaviour. In particular, the $p$-variation test diagnoses the asymptotic behaviour of 
\begin{eqnarray*}
	V_p^n(t) = \sum_{k=1}^{[n t]}|\phi(k/n) - \phi([k-1]/n)|^p.
\end{eqnarray*} 
This easily computable statistics measures the roughness of the process $\phi$, tuning into finer and finer partitions with increasing $n$. For $p=1$ the statistics reduces to the total variation and for $p=2$ it reduces to the quadratic variation. For Brownian motion where increments scale as $\sqrt{1/n}$ one obtains in the limit of $n\to \infty$ that $V_2^n(t) \sim t$, and $V_p^n(t)\to 0$ for $p>2$. Given an $\alpha$-stable process $\phi$ for some $\alpha<2$ with cumulative distribution $L_\alpha$, the statistics $V_p^n(t)$ converges for $p>\alpha$ and diverges for $p<\alpha$. It was shown that if $\phi$ is a stochastic process driven by an $\alpha$-stable process with $\alpha=p/2$ then $V_p^n(t)$ converges in distribution to $L_{1/2}$ \cite{CorcueraEtAl07,HeinImkellerPvalyukevich}. This suggests to use a Kolmogorov-Smirnov test and find the value of $p=2\alpha$ for which the empirical cumulative distribution function is closest to the target cumulative distribution function of $L_{1/2}$. To estimate the cumulative distribution function we choose to divide the time series into $15,625$ segments, each consisting of $320$ data points. The minimal Kolmogorov-Smirnov distance is then found by varying the scale parameter of the target distribution $L_{1/2}$ for each value of $p$. The value $p^\star$ for which the minimum is attained then determines $\alpha=p^\star/2$. 

Application of the $p$-variation yields $\alpha=1.28$. We have further employed different tests such as the $E(2)$-test \cite{GottwaldMelbourne16b} yielding $\alpha=1.35$, the maximum likelihood test \cite{DuMouchel73,Nolan01} yielding $\alpha=1.11$, a quantile matching test \cite{Mcculloch86} yielding $\alpha=1.73$ and linear regression of the empirical characteristic function \cite{Koutrouvelis80,Koutrouvelis81} yielding $\alpha=1.33$. For the last three methods we used the software package \cite{VeilletteMatlab}. Note that different tests typically produce different estimates for $\alpha$ for finite amounts of correlated data (see, for example, Figure 4 and 6 in \cite{GottwaldMelbourne16b}). Averaging the values for $\alpha$ obtained from these six different tests we report $\alpha\approx 1.4$. For completeness, we further estimated the tail behaviour of $d\phi$ using the method of maximum likelihood \cite{ClausetEtAl09}, and obtained $\alpha = 1.26$ for negative and for positive jumps, consistent with the results from the other tests. The skewness parameter $\beta$ was estimated to be close to zero, which means that positives jumps are equally likely as negative jumps. This was estimated by the ratio of the signs of $d\phi$.

\section*{Conclusion}
\label{sec:discussion}

We showed that the well known sLLG equation (\ref{eq:sLLG}) which is driven by Gaussian noise exhibits anomalous diffusion of the magnetisation front interface. The front interface experiences abrupt jumps leading to superdiffusion. Using the collective coordinate framework the large jumps were generated, from a purely dynamical systems perspective, by random events of unbounded drivers of the translational group variable. More concretely, the dynamical mechanism giving rise to those large jumps was identified as (near) flips of the magnetisation which induce random arbitrary large widenings of the front interface. It is pertinent to note that it was this very mechanism, that prohibited the discovery of anomalous diffusion in numerical simulations of the otherwise well studied sLLG equation: events of sufficiently large widenings require exponentially large domains to resolve the front solution for given boundary conditions. For example, given a system length $L$ for which the sLLG equation is numerically integrated, there will be an interface width $w$ such that $\tanh(w^{-1}(x-\phi))$ is not anymore compatible with any imposed Neumann boundary conditions. This is the reason why the pronounced jumps of the L\'evy process observed in Figure~\ref{fig:w} for the finite-dimensional system of SDEs for the collective coordinates are not as such visible for the solution of the full sLLG equation as depicted in Figure~\ref{fig:m1}. Whereas one could in principle subtract the drift of the interface by transforming into a frame of reference moving with the interface location, the divergence of the interface width cannot be controlled and makes it computationally impractical to resolve the jumps. The observation of rare but large jumps and the verification of the superdiffusive behaviour with sufficient statistical confidence was facilitated here by the model reduction of the collective coordinate approach.

We employed a collective coordinate approach to reduce the infinite dimensional dynamics of the SPDE to a finite set of SDEs for the collective coordinates. The choice of the collective coordinates, in particular of those describing the shape dynamics, is arbitrary. In the context of PDEs symmetry reduction can be done in a rigorous way (provided the system has a sufficient degree of hyperbolicity and supports an underlying centre manifold); see for example \cite{GolubitskyStewart}. It would be interesting to apply this machinery to the sLLG equation. We would hope that the resulting shape and group dynamics  adheres to the underlying dynamical picture of diffusion only active along the neutrally stable group dynamics unlike in the collective coordinate system (\ref{eq:CCA})-(\ref{eq:CCphi}) derived here. This more rigorous treatment of the reduction of the underlying SPDE would then be able to eliminate possible doubts that the anomalous diffusion is only of a transitory nature which currently cannot be excluded with certainty.

To further classify the underlying mechanism of the superdiffusion it would be interesting to establish to which universality class the observed anomalous diffusion belongs to by investigating the cumulant generating function associated with the collective coordinate system as introduced by \cite{StellaEtAl23a,StellaEtAl23b} who distinguish between the Richardson class \cite{Richardson26} and the Fisher class \cite{FisherSkykes59}.

Although the physical context of the sLLG equation is not the focus of this study we briefly hint to eventual physical implications of the emergent superdiffusion and discuss a possible link with the Barkhausen effect \cite{Barkhausen19}. This well known effect describes the phenomenon whereby discrete jumps in magnetization produce an audible crackling noise when a ferromagnetic material is exposed to a slowly varying external magnetic field. Such discrete jumps are known to occur in iron contaminated with impurities, where a domain wall 
can get pinned or suddenly de-pinned as it passes through these material defects. Our results suggest that such discrete jumps may not only occur due to topological constraints such as impurities but may occur generically in magnetization fronts under thermal fluctuations.


\section*{Acknowledgments}
We thank Ben Goldys and Ian Melbourne for valuable discussions. 


\section*{Appendix: Illustration of statistical limit laws}
\label{sec:Illustration}
The emergent anomalous diffusion in the sLLG equation relies on statistical limit laws. Briefly, appropriately scaled integrals (or sums) of certain time series converge to stochastic processes. In the following we illustrate this convergence to stochastic processes using entirely deterministic dynamics. For illustrative purposes we restrict to time-discrete systems. We consider systems of the form
\begin{align}
	v_{n+1} &= v_n +  \Phi(x_n)
	\label{eq:v}
	\\
	x_{n+1}&= f(x_n),
	\label{eq:x}
\end{align}
where $f(x)$ generates chaotic dynamics with an ergodic invariant probability measure $\mu$ and the driver $\Phi(x)$ is assumed to have zero mean, i.e.
\begin{align*}
	\lim_{n\to\infty} \frac{1}{n}\sum_{j=0}^{n-1} \Phi(x_j) = \int \Phi(x)d\mu = 0.
\end{align*}
We can compactly express $v_n$ as  
\begin{align*}
	v_n = \sum_{j=0}^{n-1} \Phi(x_j).
\end{align*}
The property of $\Phi(x)$ being mean-zero translates to $v$ and the question we are interested in here is how fluctuations about this mean behave. In Section~\ref{sec:BM} we discuss an example where the fluctuations are governed by the Central Limit Theorem (CLT) and in the limit $n\to \infty$ $v$ (if appropriately scaled by $1/\sqrt{n}$) exhibits Brownian motion. Depending on the observable $\Phi$ or/and the character of the driving chaotic dynamics $f(x)$ the CLT may cease to be valid and can be replaced by a modified CLT leading to superdiffusion and $v_n$ (if appropriately scaled by $1/n^\alpha$) exhibiting $\alpha$-stable dynamics. There are two known dynamical mechanisms such that $v$ may exhibit superdiffusion: either the driving $x$-dynamics is strongly chaotic with strong decay of correlations but the observable $\Phi(x)$ is not square-integrable \cite{Gouezel08} or the driving $x$-dynamics is weakly chaotic and intermittent with a sufficiently slow decay of correlations and $\Phi(x)$ may be bounded and hence square-integrable \cite{Zweimuller03,Gouezel04,TyranKaminska10a,TyranKaminska10b,MelbourneZweimueller12,ChevyrevEtAl20}. We discuss examples illustrating both of these mechanisms in Sections~\ref{sec:unbnd} and \ref{sec:intermittent}, respectively. 

In the context of the sLLG equation and symmetry groups discussed in the main text, $v$ takes the role of the group variable $\phi$ corresponding to translational symmetry, and the driving dynamics is the shape-dynamics.


\subsection{Generating Brownian motion}
\label{sec:BM}
\begin{figure}
	\centering
	\includegraphics[width=0.45\linewidth]{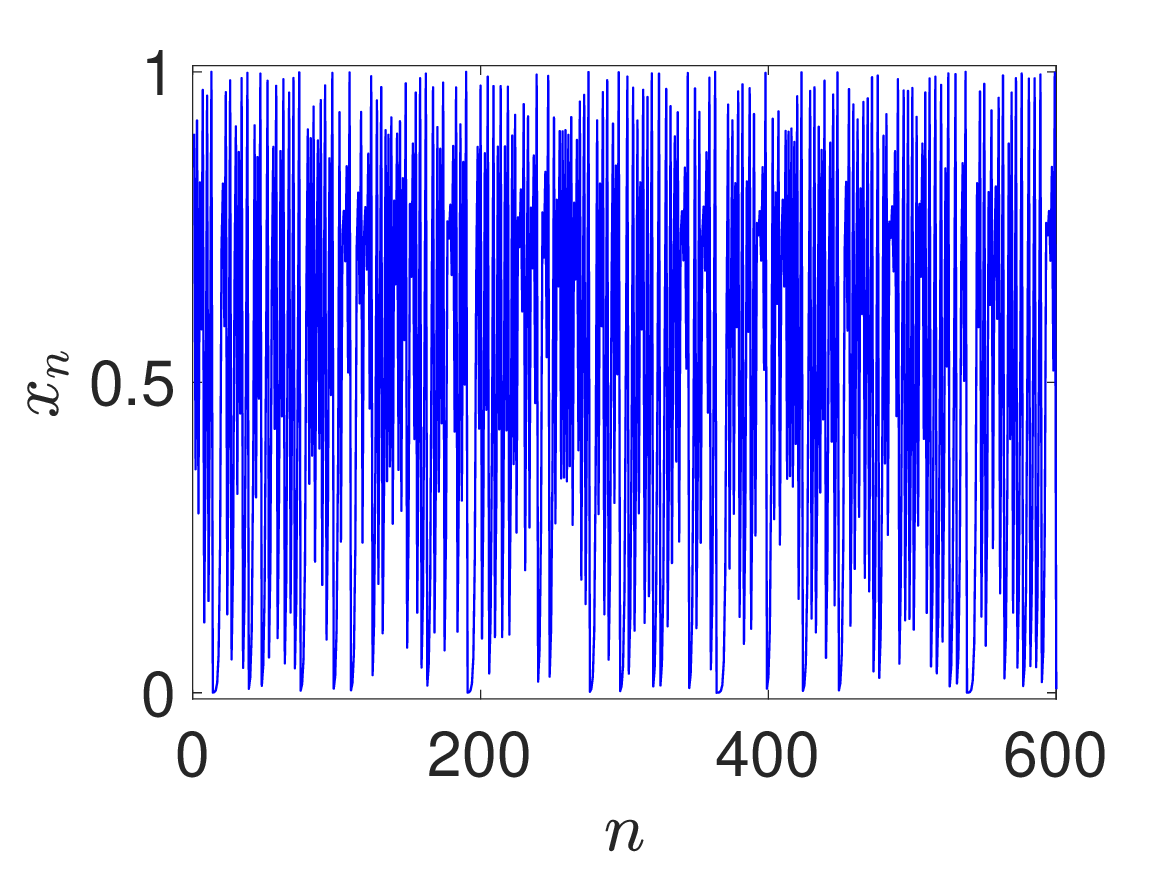}\\
	\includegraphics[width=0.45\linewidth]{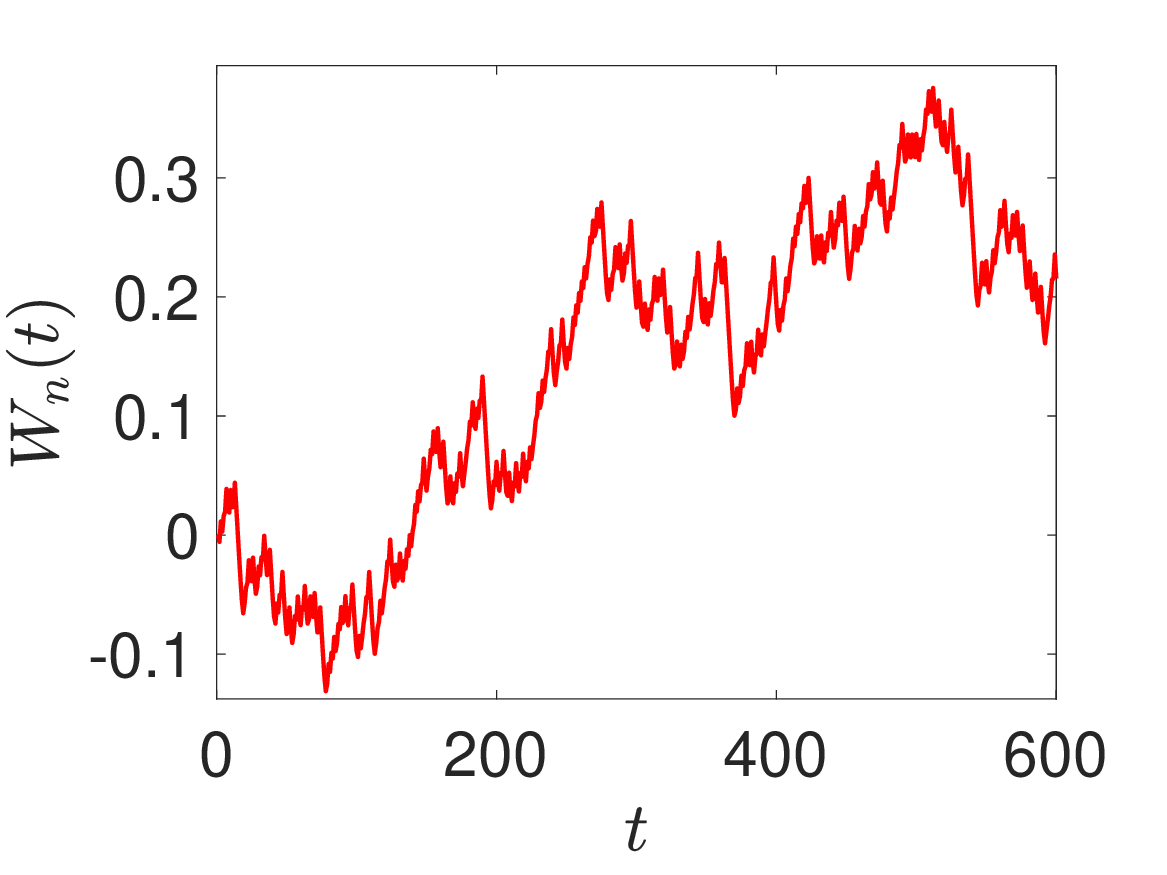}
	\caption{Simulation of (\ref{eq:logv0})--(\ref{eq:log}). Top: dynamics of the strongly chaotic logistic map. Bottom: dynamics of $W_n(t) = v_{nt}/\sqrt{n}$ with $n=600$ exhibiting asymptotically Brownian motion. To allow for a direct comparison with $x_n$ we record $t$ in units of $1/n$.}
	\label{fig:CLT}
\end{figure}

Before embarking on dynamic mechanisms to generate superdiffusion we show how to generate Brownian motion. As a particular instance of the dynamical system (\ref{eq:v})--(\ref{eq:x}) consider the following dynamical system
\begin{align}
	v_{n+1} &= v_n +  (x_n-\tfrac12)
	\label{eq:logv0}
	\\
	x_{n+1}&=4x_n(1-x_n).
	\label{eq:log}
\end{align}
The driving dynamics is the strongly chaotic logistic map (\ref{eq:log}) which enjoys exponential decay of correlations for sufficiently smooth observables and has mean $\bar x = \tfrac12$. The dynamics (\ref{eq:logv0}) for $v$ amounts to adding up a mean-zero observable leading to
\begin{align}
	\label{eq:logv}
	v_n = \sum_{j=0}^{n-1} \Phi(x_j).
\end{align}
with mean-zero $\Phi(x_j) = x_j-\tfrac12$. The statistical properties of the strongly chaotic logistic map imply a CLT with
\begin{align}
	\frac{1}{\sqrt{n}} v_n \to_d \xi,
	\label{eq:CLT}
\end{align}
where $\xi$ is a normally distributed variable $\xi\sim{\mathcal{N}}(0,\sigma^2)$ the variance of which is given by a Green-Kubo formula with 
\begin{align}
	\sigma^2=\int \Phi^2\,d\mu+2\sum_{n=1}^\infty \int \Phi(x)\;\Phi\circ f^n(x)\,d\mu.
	\label{eq:sig}
\end{align}
Note that the randomness stems from the choice of the initial condition $x_0$ which has to be randomly drawn. Moreover, besides the convergence of $v_n$ in law to a Gaussian random variable, the continuous function defined by
\begin{align}
	W_n(t) = \frac{1}{\sqrt{n}} v_{nt} = \frac{1}{\sqrt{n}} \sum_{j=0}^{nt-1} \Phi(x_j),
	\label{eq:W}
\end{align} 
for discrete $t=0,\tfrac1n,\tfrac2n,\cdots$, and where we linearly interpolate to obtain a continuous function converges weakly to a Wiener process $W$ with variance $\sigma^2 t$.  

To illustrate these heuristic arguments we show in Figure~\ref{fig:CLT} a time series obtained from simulating the system (\ref{eq:logv0})--(\ref{eq:log}), where it becomes visible how by integrating strongly chaotic variables, $W_n$ approaches Brownian motion.

For rigorous proofs of these statements and observations the reader is referred to \cite{MelbourneNicol09,MelbourneNicol05,Gouezel10,GottwaldMelbourne13c}. We remark that these works made rigorous the ingenious idea, that the climate can be treated as a stochastic dynamical system as a result of experiencing the integrated effect of fast moving weather systems, for which Klaus Hasselmann received the 2022 Nobel Prize in physics.


\subsection{Generating $\alpha$-stable processes}

For the variance (\ref{eq:sig}) of the Brownian motion to be defined the auto-correlation function needs to be summable. This can be violated either by choosing observables which are not square-integrable (leading to divergence of the first integral in the expression (\ref{eq:sig}) for $\sigma^2$) or by driving dynamics which allows for correlations over arbitrary large time intervals (leading to the divergence of the sum in the expression (\ref{eq:sig})). 


\subsubsection{Unbounded observables}
\label{sec:unbnd}
We consider again a system driven by the strongly chaotic logistic map, but now the observable $\Phi(x)$ is unbounded and not square-integrable. Concretely, we consider
\begin{align}
	v_{n+1} &= v_n +  (x_n^{-\gamma}-\bar x)
	\label{eq:logunbndv} \\
	x_{n+1}&=4x_n(1-x_n),
	\label{eq:logunbnd}
\end{align}
where $\gamma<\tfrac12$ to render $\Phi(x)=x_n^{-\gamma}-\bar x$ non-square-integrable and $\bar x=\bar x(\gamma)=\Gamma(\tfrac12-\gamma)/(\sqrt{\pi}\Gamma(1-\gamma))$ is chosen to ensure that $\Phi(x)$ is mean zero. The variable $v$ now experiences strong kicks whenever the logistic map evolves near $x_n=0$. These strong kicks lead to quasi-discrete jumps in $v$. As for convergence to Brownian motion, an appropriate scaling of $v_n$ is required to obtain convergence to an $\alpha$-stable process for $n \to \infty$; in particular 
\begin{align}
	L_n(t) =\frac{1}{n^{\frac{1}{\alpha}}} v_{nt} =  \frac{1}{n^{\frac{1}{\alpha}}} \sum_{j=0}^{nt-1} \Phi(x_j)
	\label{eq:L}
\end{align}
asymptotically approximates an $\alpha$-stable process where $\alpha = \alpha(\gamma)$. In our example the $\alpha$-stable process is asymmetric with $\beta=1$ allowing only for positive jumps. 

\begin{figure}
	\centering
	\includegraphics[width=0.34\linewidth]{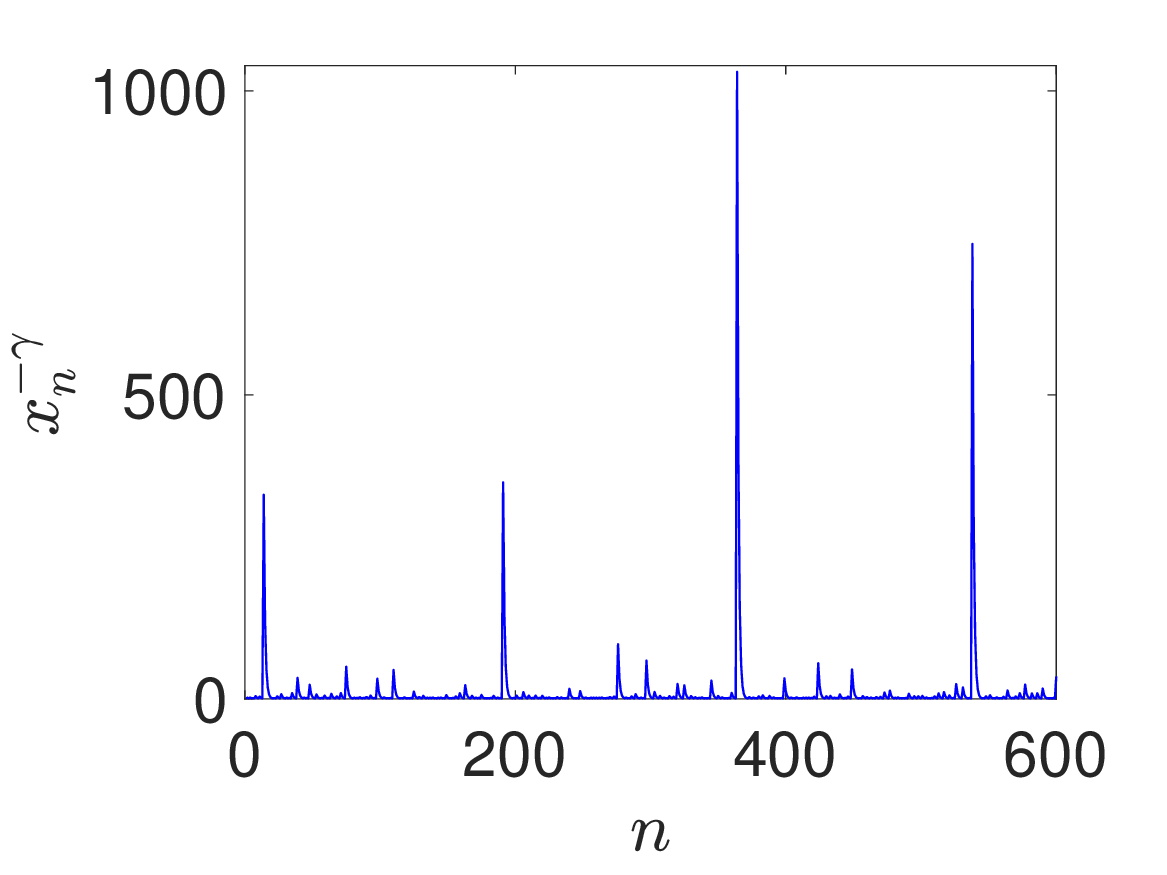}\\
	\includegraphics[width=0.34\linewidth]{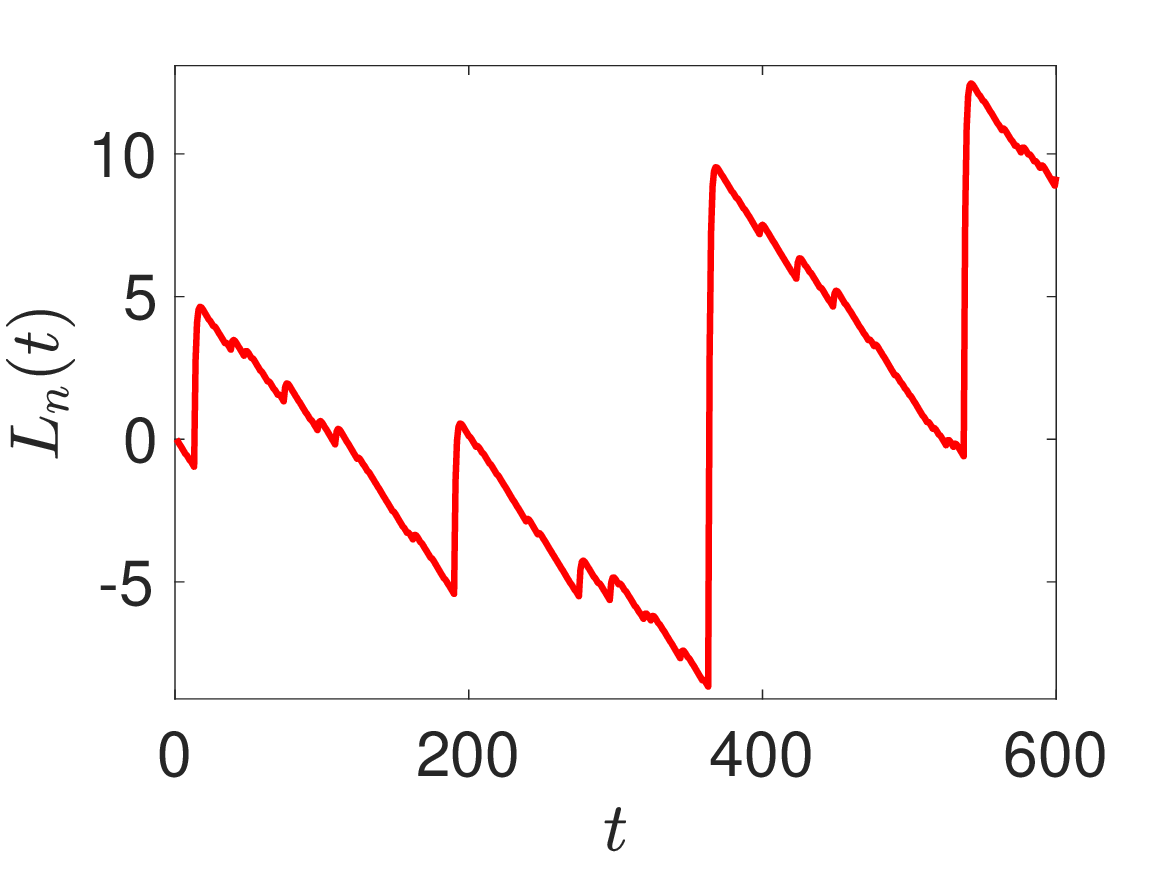}
	\caption{Simulation of (\ref{eq:logunbndv})--(\ref{eq:logunbnd}) with $\gamma = 0.4$. Top: dynamics of the strongly chaotic logistic map for the observable $x^{-\gamma}$ experiencing spikes whenever $x_n\ll 1$. Bottom: dynamics of $L_n(t)  =v_{nt}/n^{1/\alpha}$ with $n=600$ and $\alpha=1/2\gamma=1.25$ exhibiting asymptotically superdiffusion. To allow for a direct comparison with $x_n$ we record $t$ in units of $1/n$.}
	\label{fig:Levy_unbnd}
\end{figure}

This is illustrated in Figure~\ref{fig:Levy_unbnd} where we see clearly how spikes in $x_n^{-\gamma}$ generate jumps in $v_n$ (and hence in $L_n(t)$), the size of which depend on the magnitude of the driving spike. 

For completeness we determine the value of $\alpha$ corresponding to a given exponent $\gamma$. It is well-known that a mean-zero random variable $Z^{-\gamma}$ lies in the domain of attraction of an $\alpha$-stable law, i.e. that $\sum_{j=0}^{n-1} Z_j^{-\gamma}/n^\alpha$ converges in law to an $\alpha$-stable random variable, if it has sufficiently strong heavy tails \cite{GnedenkoKolmogorov}. Applying this to our case, we have convergence when $P(Z^{-\gamma}>z)=P(Z<z^{-\tfrac1\gamma})\sim z^{-\alpha}$ with $\alpha\in(0,2)$. Given the probability density function of the logistic map $p(x) = 1/\pi\sqrt{x(1-x)}$ we readily evaluate
\begin{align}
	P(x^{-\gamma}>\zeta)=P(x<\zeta^{-\tfrac1\gamma})=\frac{2}{\pi}\sin^{-1} \sqrt{\zeta^{-\tfrac1\gamma}} \sim \zeta^{-\frac{1}{2\gamma}}
\end{align} 
for large values of $\zeta$, and hence we obtain $\alpha=1/(2\gamma)$. 

Rigorous proofs of these statements and observations exist for uniformly expanding maps and the reader is referred to \cite{Gouezel08}. The stochastic counterpart to this deterministic mechanism is given by correlated additive and multiplicative (CAM) noise, which is often used in the climate sciences to model unresolved atmospheric dynamics. CAM noise exhibits intermittent unbounded spikes and was shown rigorously to lie in the domain of attraction of $\alpha$-stable noise \cite{KuskeKeller01}.

It is this mechanism of unbounded observables which we identified to be active in the sLLG equation; recall that the dynamics of the translational variable $\phi$ is driven by the interface width of the magnetisation front $w$ which can become arbitrarily large (cf Equation (10) in the main text). 


\subsubsection{Intermittent dynamics}
\label{sec:intermittent}
\begin{figure}[thb]
	\centering
	\includegraphics[width=0.32\linewidth]{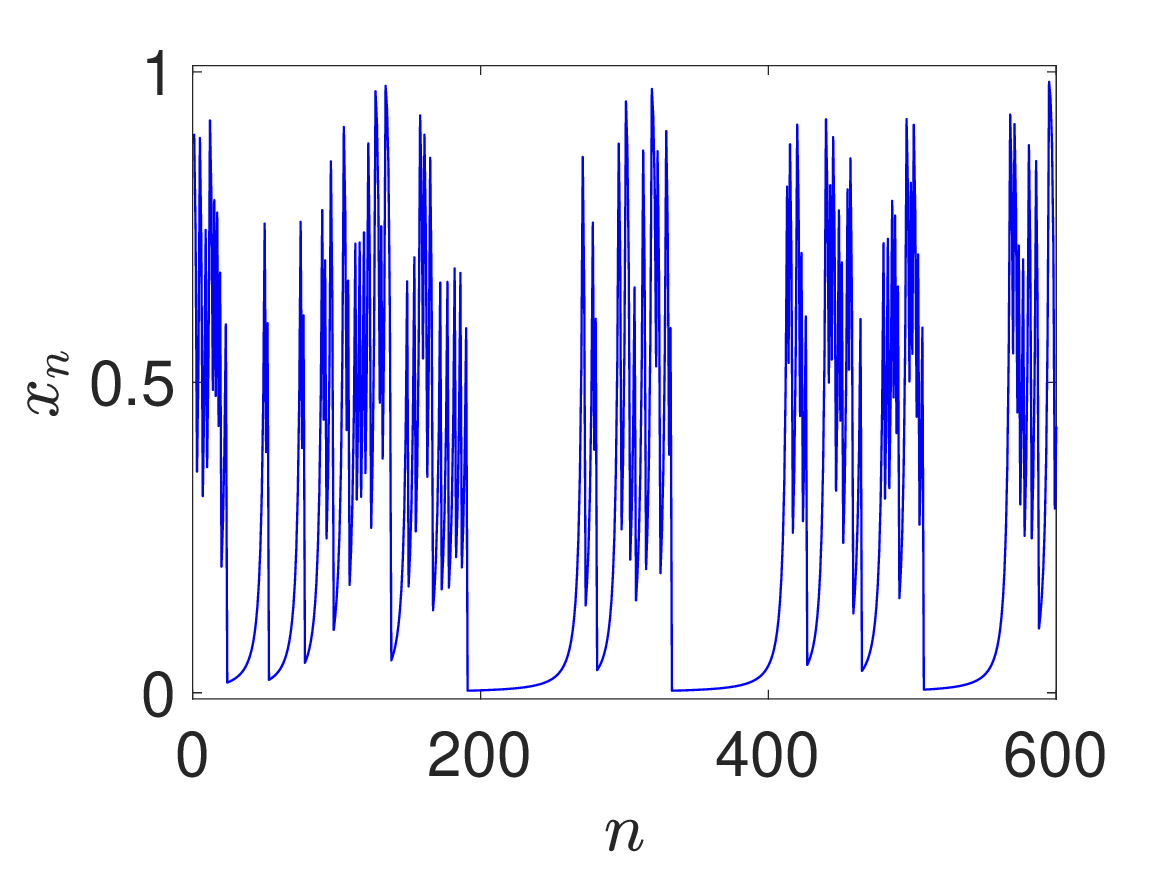}\\
	\includegraphics[width=0.32\linewidth]{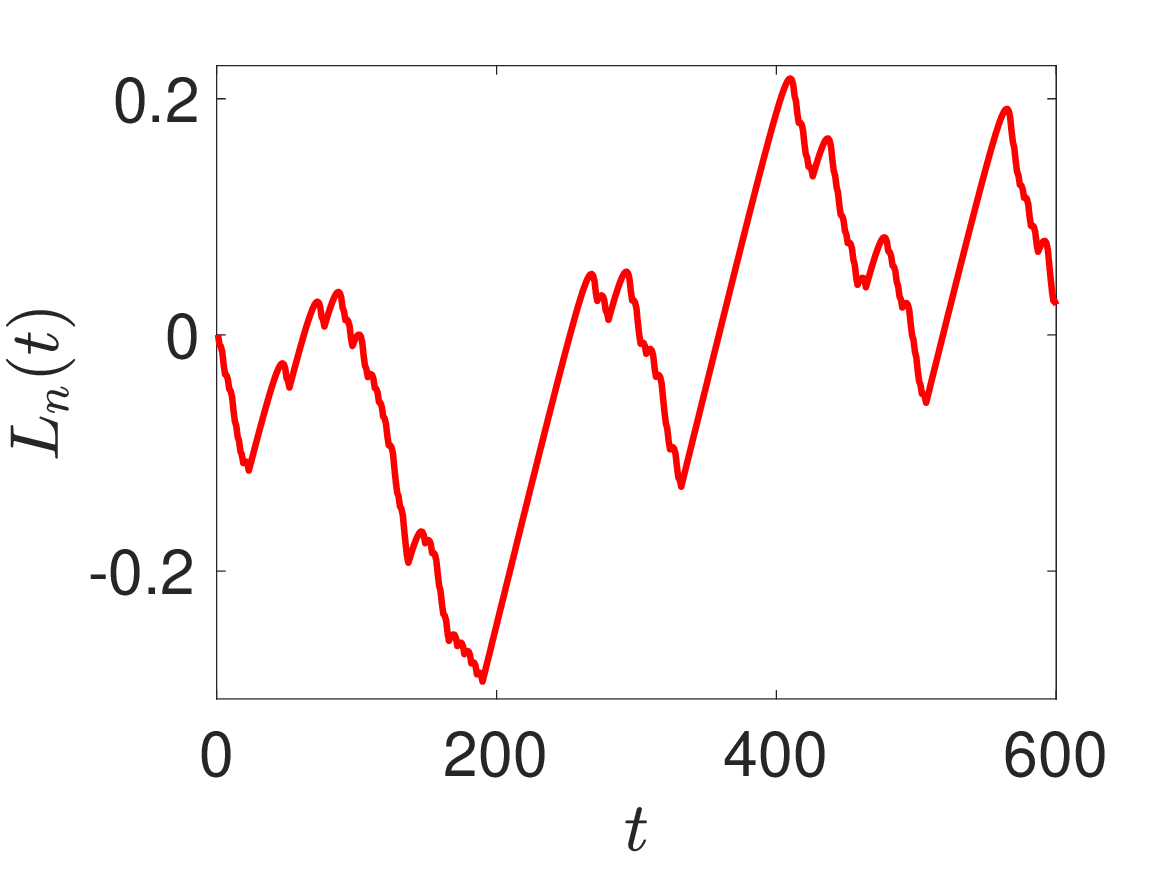}\\
	\includegraphics[width=0.32\linewidth]{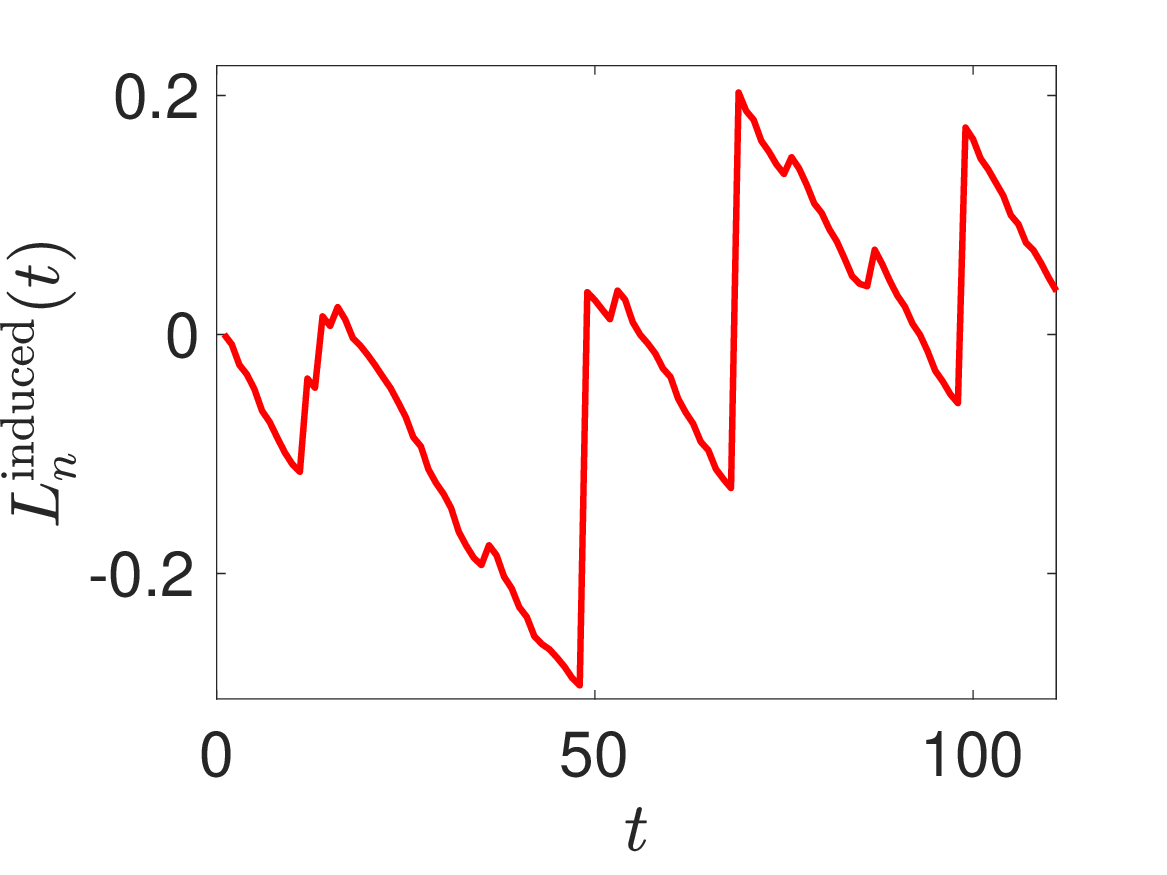}
	\caption{Simulation of (\ref{eq:intermittv})--(\ref{eq:intermittv}) with $\gamma = 0.7$ and $x^\star = 0.5$. Top: dynamics of the weakly chaotic intermittent map. Middle: dynamics of $L_n(t)  = v_{nt}/n^{1/\alpha}$ with $n=600$ and $\alpha=1/\gamma=1.43$ exhibiting asymptotically superdiffusion. To allow for a direct comparison with $x_n$ we record $t$ in units of $1/n$. Bottom: dynamics of the induced variable $L_n^{\rm{induced}}(t)$ obtained by stopping the time whenever the driving $x$-dynamics enters a laminar period with $x<0.5$.} 
	\label{fig:Levy_intermittent}
\end{figure}

For completeness we illustrate how superdiffusion can be generated by intermittent dynamics in which the dynamics exhibits prolonged laminar periods interspersed with intermittent chaotic bursts. For intermittent behaviour where the laminar periods become too long for the correlation function to be summable, the CLT breaks down \cite{GaspardWang88}. We now show how this may give rise to $\alpha$-stable processes and superdiffusion. We consider
\begin{align}
	v_{n+1} &= v_n - (x_n-\bar x)
	\label{eq:intermittv}
	\\
	x_{n+1}&=\begin{cases}
		x_n(1+2^\gamma x_n^\gamma) & 0\le x_n\le \frac12 \\ 
		2x_n-1 & \frac12\le x_n\le1 
	\end{cases}
	\label{eq:intermitt}
\end{align}
with $\gamma\in (\tfrac12,1)$ and where $\bar x=\bar x(\gamma)$ is chosen to ensure that $\Phi(x)$ has zero mean. The driving intermittent dynamics (\ref{eq:intermitt}) is a chaotic Pomeau-Manneville map \cite{PomeauManneville80,LiveraniEtAl99} which has a neutrally stable fixed point at $x=0$. When the dynamics evolves near this neutrally stable fixed point it experiences sticky laminar behaviour with $x\approx 0$ for some time. For $\gamma\ge\frac12$, the stickiness is strong enough to render the autocorrelation function nonsummable \cite{Hu04}. For $\gamma> \tfrac12$ the CLT is replaced by a modified CLT where 
$n^{-\gamma}v_n$ converges in law to an $\alpha$-stable random variable with $\alpha=1/\gamma$ \cite{Gouezel04}, and the deterministic dynamics of the appropriately scaled $v_n$, i.e. the process $L_n(t)$ defined above in (\ref{eq:L}), converges for $n\to\infty$ to a corresponding $\alpha$-stable process, exhibiting superdiffusion \footnote{For $\gamma=\tfrac12$ a different scaling with $1/\sqrt{n \log n}$ is required instead to achieve convergence.} .

The heuristic mechanism how intermittent dynamics can induce jumps in $v_n$ when summed up is nicely illustrated in the numerical simulations of the system (\ref{eq:intermittv})--(\ref{eq:intermitt}), depicted in Figure~\ref{fig:Levy_intermittent}. Whenever the driving $x$-dynamics becomes laminar with $x\approx 0$ we have $\Phi(x)=-(x-\bar x) \approx \bar x$ (see top row). Hence provided $\Phi(x=0)=\bar x$ is non-zero at the neutral fixed point, the $v$ dynamics evolves ballistically with constant speed $\bar x$ according to $v_{n+1}\approx v_n + \bar x$ and during the laminar periods we can approximate $v_n\sim \bar x n$. These ballistic flights can become arbitrarily large provided the laminar periods can become arbitrarily long - this is possible with a non-neglibile probability for $\gamma>\tfrac12$. The ballistic dynamics of $v$ is clearly seen in the middle row of Figure~\ref{fig:Levy_intermittent} where we show the scaled process $L_n(t)$. The ballistic flights correspond in the long time limit to the discrete jumps of an $\alpha$-stable process (with $\alpha=1/\gamma$), each jump consisting of arbitrary many small steps of size $\bar x$ each. To obtain a more distinct process one can employ what is called {\em{inducing}} which amounts to stopping time whenever the dynamics is in the laminar region with $x\le 0.5$ and only record $v(x)$ when the driving dynamics is in the expanding strongly chaotic region with $x\ge 0.5$ (bottom row). After inducing all the small steps during the ballistic flight are being ignored and we obtain actually discrete jumps. Both processes, the uninduced (middle row) and the induced one (bottom row) asymptotically converge to $\alpha$-stable processes, however in different topologies (see \cite{ChevyrevEtAl20} for precise statements); from a non-mathematical point of view it is the degree of squinting required to make $v_n/n^\gamma$ look like an $\alpha$-stable process with discrete jumps which distinguishes the two cases.

For rigorous proofs of these statements and observations, as well as for explicit formulae on how to calculate the parameters of the limiting $\alpha$-stable process, the reader is referred to \cite{Zweimuller03,Gouezel04,TyranKaminska10a,TyranKaminska10b,MelbourneZweimueller12,GottwaldMelbourne13c,ChevyrevEtAl20}.



\end{document}